\documentclass{article}

\usepackage{arxiv}

\usepackage[utf8]{inputenc} 
\usepackage[T1]{fontenc}    
\usepackage{hyperref}       

\usepackage{url}            
\usepackage{booktabs}       
\usepackage{amsfonts}       
\usepackage{amsmath}
\usepackage{amssymb}
\usepackage{nicefrac}       
\usepackage{microtype}      
\usepackage{graphicx}
\usepackage{multirow}
\usepackage{adjustbox}     
\usepackage[section]{placeins}  
\usepackage[ruled,linesnumbered]{algorithm2e}

\SetArgSty{textnormal}   
\graphicspath{ {./images/} }

\newcommand{\supplement}{the supplementary material (ancillary file \texttt{supplement.pdf})}

\title{Tree-based formulation for the multi-commodity flow problem}

\author{
 Simon Spoorendonk \\
  Flowty ApS\\
  Denmark\\
  \texttt{simon@flowty.ai} \\
   \And
 Bjørn Petersen \\
  Flowty ApS\\
  Denmark\\
  \texttt{bjorn@flowty.ai} \\
}

\begin{document}
\maketitle
\begin{abstract}
We revisit the tree-based formulation of the minimum-cost multi-commodity flow problem, due to
Jones et al.~\cite{jones1993multicommodity}, who found that path-based decomposition converges in
fewer master iterations and reported lower CPU times for it. The formulation represents the flow
out of each source as a convex combination of shortest-path trees, so the master problem has one
demand constraint per source $|S|$ rather than per commodity $|K|$. We re-examine it on 44
instances with up to 3.3 million commodities, three orders of magnitude beyond the scale
available to Jones et al., under five linear programming backends spanning four barrier codes,
open source and commercial, CPU and GPU. Their convergence result is confirmed: the tree-based
formulation still requires two to three times as many iterations. Their wall-clock conclusion,
however, is reversed under every backend in the regime $|S| \ll |K|$: tree-based column
generation is 1.4
to 1.9 times faster on the shifted geometric mean, reaches a factor of 99 on the instance with the
most commodities, and solves 43 or 44 of the 44 instances under every backend, where the
path-based formulation solves 38 to 42, its failures concentrated on the largest transportation
instances. Both decompositions are 15 to 28 times faster than
solving the compact model directly (5 to 25 on the instances that model solves). The
measurements identify the mechanism: the master problem accounts for 87\% to 99\% of the
runtime, and the tree-based master is up to 32 times smaller at termination. The advantage
disappears when $|S|$ approaches $|V|$, where a tree column has non-zeros in a large fraction of
the master's rows, as on the planar2500 instance. An open-source
C++ implementation accompanies the paper.

\end{abstract}

\keywords{Multi-commodity flow problem \and Column generation}

\section{Introduction}

The minimum-cost multi-commodity flow (MCF) problem is a fundamental network optimization problem 
with applications in logistics, passenger transportation, telecommunication, and energy systems 
\cite{ahuja1993network,salimifard2022mcnf}. It routes multiple origin–destination demands through 
a capacitated network at minimum cost. The MCF problem is a linear program, but the number of 
variables and constraints grows rapidly, making large-scale instances computationally challenging.

The classical edge-based formulation introduces $O(|K||E|)$ variables and $O(|K||V| + |E|)$ constraints. 
The model becomes impractically large when either the commodity set $K$ or the network grows. 
Source-based formulations aggregate by source nodes, reducing the model to $O(|S||E|)$ variables 
and $O(|S||V| + |E|)$ constraints. When $|S| \ll |K|$ this yields savings, but commodity-level 
flows must then be reconstructed \cite{ahuja1993network}. 
Path-based formulations, obtained via Dantzig–Wolfe decomposition \cite{dantzig1960decomposition}, replace edge variables by
path variables. Each commodity $k \in K$ has an exponential path set $\mathcal{P}_k$, and column 
generation dynamically adds paths \cite{desrosiers2024branch,uchoa2024column}. 
However, the master problem still contains $O(|K|)$ demand constraints, and restricted master 
problems grow prohibitively large when $|K|$ reaches millions.

Research has therefore focused on decomposition and relaxation
approaches. Stabilized column generation methods
\cite{gondzio2016primaldual}, analytic-centre cutting-plane algorithms
with active sets, and Lagrangian relaxation
\cite{babonneau2006active} have extended the tractable instance size.
Recently, Zhang and Boyd \cite{zhang2025gpu} proposed a GPU-compatible
primal-dual hybrid gradient method for an all-pair utility-maximizing variant
of the MCF problem (their objective is a strictly concave utility rather than a linear
cost), itself based on a destination-aggregated formulation, achieving
orders-of-magnitude speed-ups over interior-point solvers; whether this carries over to the
min-cost MCF problem remains to be investigated.
Destination-aggregated formulations have also been investigated in
communication and optical networks \cite{yin2020unified}, and origin-aggregated
formulations are standard in traffic assignment \cite{bargera2002origin}. In each
case, scalability comes from exploiting problem structure rather than introducing flow
variables that are not needed.

To address the master problem bottleneck, we revisit the tree-based
decomposition of the MCF problem, first studied by Jones et al.~\cite{jones1993multicommodity}.
Their formulation derives from the source-based formulation by representing flows as combinations
of trees rooted at source nodes, so the master problem has $O(|S|)$ demand constraints and
$O(|E|)$ capacity constraints, against $O(|K|)$ demand constraints for the path-based
formulation. Jones et al.\ found that path-based decomposition converged
faster at the scales they tested and concluded paths were preferable, noting (p.96) that
``larger and sparser master problems do not necessarily entail longer solution times at each
iteration.'' We re-examine this at modern scale, where $|K|$ reaches millions and
$|S| \ll |K|$, and find that the trade-off reverses: the compactness of the $O(|S|)$ tree
master problem outweighs the per-iteration cost advantage of the path-based formulation. The
tree-based master is also the smaller of the two in columns at termination, by a factor of 32 on
the transportation family, whereas Jones et al.\ report the opposite on two of their three
largest instances.
To our knowledge, we are the first to demonstrate this
computationally at million-commodity scale: in the regime $|S| \ll |K|$, the source-based
tree decomposition solves instances substantially faster than the path-based formulation, and
solves instances the path-based formulation cannot finish.

The tree-based formulation itself is due to Jones et al.; the contribution of this paper is the
computational finding at scale, together with the techniques of \autoref{sec:colgen} that make
the formulation solvable at this size.

The contributions of this paper are:
\begin{itemize}
  \item \textbf{Empirical reversal at scale.} On instances with up to 3.3 million commodities,
  three orders of magnitude larger in $|K|$ than those Jones et al.\ could test, the ordering
  they reported is reversed: in the regime $|S| \ll |K|$ the source-based tree decomposition
  outperforms path-based decomposition on the shifted geometric mean under every one of the
  five LP backends. The reversal is an aggregate one (under two backends the path-based
  formulation is still faster on a majority of individual instances, but the margins by which it
  wins are small and those by which it loses are large), and it is not attributable to the degree of source aggregation but to
  absolute scale, as \autoref{sec:conclusion} sets out. We characterize the regimes empirically
  via a heatmap of speed-up as a function of $|K|$ and $|S|/|K|$.
  \item \textbf{Solver-independent evidence.} We repeat the entire study under five linear
  programming backends spanning open-source and commercial codes and both CPU and GPU execution.
  All five reproduce the tree/path ordering, which separates the formulation effect from the
  characteristics of any single barrier implementation, a control that computational
  comparisons of decomposition schemes rarely apply.
  \item \textbf{Pricing acceleration techniques.} We introduce a bounded stopping criterion for
  shared-source pricing, which cuts a single multi-target shortest-path run short on the maximum
  of the demand duals still outstanding, together with a demand-weighted analogue for tree
  pricing that is not a rescaling of it. Both come with correctness arguments, and
  \autoref{sec:ablation} measures the criterion in a paired ablation over both the instance family
  and the LP backend. We also extend the A*-based pricing of Lienkamp and
  Schiffer~\cite{lienkamp2024column} to the multi-target tree setting using a precomputed reverse
  multi-source Dijkstra heuristic with $O(|V|)$ storage.
  \item \textbf{Master/pricing balancing.} We propose a two-mode strategy, master-problem-heavy
  or pricing-heavy, selected according to which component is the computational bottleneck on an
  instance family, incorporating lazy capacity constraints
  \cite{gondzio2016primaldual,babonneau2006active} and the pricing filter of Lienkamp and
  Schiffer~\cite{lienkamp2024column}.
\end{itemize}

An open-source C++ implementation, together with the instance-conversion scripts and the scripts
that rebuild every table and figure from the committed result files, is available at
\url{https://github.com/spoorendonk/mcfcg}. Every number reported here was produced by release
v0.1.0, which is archived under its own version DOI~\cite{mcfcg}; the repository is a moving
reference and the archive is what reproduces these results.

The remainder of the paper is organized as follows.
\autoref{sec:problem} introduces the formulations of the
MCF problem including the tree-based formulation of Jones et al.~\cite{jones1993multicommodity}. \autoref{sec:colgen}
describes the column generation algorithm used for the path- and tree-based formulations.
\autoref{sec:experiments} presents computational experiments on
large-scale instances, and \autoref{sec:conclusion} concludes with
directions for future research.

\section{Problem statement}
\label{sec:problem}

We present four formulations: the edge-based formulation, the source-based formulation, the path-based
formulation commonly solved via column generation, and the tree-based formulation \cite{jones1993multicommodity}.

We consider a directed network $G=(V,E)$ with node set $V$ and edge set
$E$. Each edge $e \in E$ has capacity $u_e \geq 0$. The set of
commodities is $K$, with each commodity $k \in K$ defined by source
$s_k \in V$, sink $t_k \in V$, and demand $d_k > 0$. Sending one unit
of any commodity over edge $e$ incurs a cost $c_e \geq 0$. The objective of the MCF
problem is to route all demands at
minimum cost without exceeding edge capacities.

\subsection{Edge-based formulation}
The classical edge-based linear program introduces variables
$f^k_e$ for the flow of commodity $k$ on edge $e$. The model is

\begin{align}
  \min \quad & \sum_{k \in K} \sum_{e \in E} c_e f^k_e \label{eq:edge-obj} \\
  \text{s.t.} \quad
  & \sum_{e \in \delta^+(i)} f^k_e - \sum_{e \in \delta^-(i)} f^k_e =
  \begin{cases}
    d_k, & i = s_k, \\
    -d_k, & i = t_k, \\
    0, & \text{otherwise},
  \end{cases}
  && \forall i \in V, \; k \in K, \label{eq:edge-balance} \\
  & \sum_{k \in K} f^k_e \leq u_e, && \forall e \in E, \\
  & f^k_e \ge 0, && \forall e \in E,\; \forall k \in K, \label{eq:edge}
\end{align}
where $\delta^+(i)$ and $\delta^-(i)$ denote the outgoing and incoming
edges of node $i$, respectively. This formulation has $O(|K||E|)$
variables and $O(|K||V| + |E|)$ constraints.

\subsection{Source-based formulation}
Source-based formulations aggregate flows by source node. Let
$S = \{s_k : k \in K\}$ denote the set of unique source nodes, let
$K_s = \{k \in K : s_k = s\}$ denote the commodities with source $s \in S$,
and let $D_s = \sum_{k \in K_s} d_k$ denote the total demand originating at $s$.
For each source $s \in S$, let $f^s_e$ denote the total flow on edge $e$
originating at $s$. The formulation is
\begin{align}
  \min \quad & \sum_{s \in S} \sum_{e \in E} c_e f^s_e 
  \label{eq:source-obj} \\
  \text{s.t.} \quad
  & \sum_{e \in \delta^+(i)} f^s_e - \sum_{e \in \delta^-(i)} f^s_e = d^s_i,
  && \forall i \in V,\; \forall s \in S, \\
  & \sum_{s \in S} f^s_e \leq u_e,
  && \forall e \in E,  \\
  & f^s_e \ge 0, && \forall e \in E,\; \forall s \in S, \label{eq:source}
\end{align}
where the aggregate demand of source $s$ at node $i$ is
$d^s_i = D_s$ if $i = s$ and $d^s_i = -\sum_{k \in K_s : t_k = i} d_k$ otherwise;
the sum is empty, and $d^s_i$ zero, at nodes where no commodity from $s$ has its sink.

This model has $O(|S||E|)$ variables and $O(|S||V| + |E|)$ constraints.
When $|S| \ll |K|$, this reduces the model size by a factor of about $|K|/|S|$ relative to the edge-based formulation.
The two models have the same optimal value: arc costs do not depend on the commodity, so any
source-aggregated flow decomposes into commodity flows of equal cost, and conversely summing a
commodity solution over each source gives a source-aggregated one. This is what allows the
tree-based and path-based masters, which are Dantzig--Wolfe reformulations of the source-based and
edge-based models respectively, to be compared against a single reference objective.
The trade-off is that the commodity-level flow decomposition is lost and must be reconstructed
in a post-processing step: standard flow decomposition writes the flow $f^s$ out of each source
as at most $|E|$ paths and cycles in $O(|V||E|)$ time \cite{ahuja1993network}, and the cycles,
which carry no demand, are discarded, which cannot raise the cost since $c \ge 0$.

\subsection{Path-based formulation}
Path-based formulations rely on a Dantzig--Wolfe decomposition of the edge-based formulation:
each commodity $k$ defines a sub-problem whose feasible region is the set of all $s_k$--$t_k$
flows of value $d_k$, and the extreme points of that polytope are the simple $s_k$--$t_k$
paths \cite{desrosiers2024branch,uchoa2024column,barnhart1994column,barnhart2000branch}. For each
commodity $k \in K$, let $\mathcal{P}_k$ be the set of all such paths, and introduce variables
$x^k_p \geq 0$ denoting the flow of commodity $k$ on path $p \in \mathcal{P}_k$. The model is
\begin{align}
  \min \quad & \sum_{k \in K} \sum_{p \in \mathcal{P}_k} c_p x^k_p 
  \label{eq:path-obj} \\
  \text{s.t.} \quad
  & \sum_{p \in \mathcal{P}_k} x^k_p = d_k, && \forall k \in K, \label{eq:path-balance} \\
  & \sum_{k \in K} \sum_{p \in \mathcal{P}_k : e \in p} x^k_p \leq u_e,
  && \forall e \in E, \label{eq:path-cap}
   \\
  & x^k_p \geq 0, && \forall k \in K, \; p \in \mathcal{P}_k, \label{eq:path}
\end{align}
where $c_p = \sum_{e \in p} c_e$ is the cost of path $p$. The formulation has $O(|K|+|E|)$
constraints; the path set $\mathcal{P}_k$ is exponential in $|V|$, and the model is solved by
column generation, see \autoref{sec:colgen}.

\subsection{Tree-based formulation}
\label{sec:tree}
Jones et al.~\cite{jones1993multicommodity} study a destination-specific decomposition (DSP) in
which a commodity aggregates all flows arriving at a common destination, and the subproblem is a
shortest-path tree converging into that destination. We study the origin-based counterpart, in
which a commodity aggregates all flows leaving a common source and the subproblem is a
shortest-path tree rooted at it. Under the standard assumption that arc costs are not origin- or
destination-specific (which holds for all instances in our experiments), the two are equivalent
by reversing the direction of aggregation, and Jones et al.\ themselves apply their formulation
in this origin-flow direction on one of their data sets~\cite[p.~111]{jones1993multicommodity}.
In practice one aggregates on the smaller of the two index sets; only the source-based variant is
tested here. The tree-based formulation is the Dantzig--Wolfe decomposition of the source-based
formulation, representing the flow out of each source as a convex combination of
\emph{shortest-path trees} rooted at it. The subproblem polytope is unbounded (its recession cone
is the cycle cone), so the reformulation carries extreme rays in general. We omit them: costs are
non-negative, so every cycle both raises the objective and consumes capacity, and any feasible
solution has a cycle-free counterpart of no greater cost. The same remark applies to the
path-based master.

For each source $s \in S$, let $\mathcal{T}_s$ denote the set of
trees rooted at $s$ that contain all sinks $\{t_k : k \in K_s\}$ of the
commodities $K_s$ with source $s$. A decision variable $x^s_\tau \geq 0$ is
associated with each tree $\tau \in \mathcal{T}_s$, indicating the fraction
of flow routed through $\tau$. Let $e \in \tau$ be an edge in the tree $\tau$ and let $p_k$ be the unique path from $s$ to $t_k$ in the tree $\tau$ for $k \in K_s$, where we suppress the dependence of $p_k$ on $\tau$ for brevity.
The model is

\begin{align}
  \min \quad & \sum_{s \in S} \sum_{\tau \in \mathcal{T}_s} c_\tau x^s_\tau
  \label{eq:tree-obj} \\
  \text{s.t.} \quad
  & \sum_{\tau \in \mathcal{T}_s} x^s_\tau = 1, 
  && \forall s \in S, \label{eq:tree-balance} \\
  & \sum_{s \in S} \sum_{\tau \in \mathcal{T}_s:\, e \in \tau} \bar{f}^e_\tau x^s_\tau \leq u_e, 
  && \forall e \in E, \label{eq:tree-cap} \\
  & x^s_\tau \geq 0, && \forall s \in S,\; \tau \in \mathcal{T}_s. \label{eq:tree}
\end{align}

Here $\bar{f}^e_\tau = \sum_{k \in K_s : e \in p_k} d_k$ is the total flow on edge $e$ in tree $\tau$,
and $c_\tau = \sum_{e \in \tau} \bar{f}^e_\tau c_e$ denotes the cost of tree $\tau$. The formulation has $O(|S|+|E|)$ constraints. The tree set $\mathcal{T}_s$ is exponential in size.

\section{Column Generation}
\label{sec:colgen}
Column generation solves a linear program whose variables are too numerous to enumerate by
alternating between a restricted master problem (RMP) over a subset of them and pricing
sub-problems that generate variables of negative reduced cost from the RMP's duals; when none
exists, the RMP solution is optimal for the full problem. The idea dates back to the
multi-commodity flow computations of Ford and Fulkerson~\cite{ford1958suggested} and has since
been applied to many problem classes \cite{desrosiers2024branch,uchoa2024column}. For the MCF
problem the path-based and tree-based formulations differ only in the pricing problem, shortest
paths against shortest-path trees. We describe each, then the implementation techniques that
exploit the network structure.

\subsection{Path generation}
For the path-based formulation the column generation procedure alternates between:

\begin{itemize}
  \item \emph{Restricted master problem (RMP):} solve
    \eqref{eq:path-obj}--\eqref{eq:path} using only a subset of paths from $\cup_{k \in K} \mathcal{P}_k$.
  \item \emph{Pricing sub-problem:} for each commodity $k$, compute a
    shortest $s_k$--$t_k$ path with respect to the reduced costs given
    by the dual variables of the RMP. If a path of negative reduced cost
    is found, it is added as a new column.
\end{itemize}

Let $\pi_k$ denote the dual variable associated with the demand
constraint of commodity $k$ \eqref{eq:path-balance}, and let $\mu_e$ denote the dual variable
associated with the capacity constraint on edge $e$ \eqref{eq:path-cap}. Since
\eqref{eq:path-cap} is a less-than-or-equal constraint in a minimization problem, the duals
satisfy $\mu_e \leq 0$. Then the reduced
cost of a path $p \in \mathcal{P}_k$ is
\begin{align}
  \bar{c}_p 
  = \sum_{e \in p} \big(c_e - \mu_e\big) - \pi_k.
  \label{eq:path-reduced-cost}
\end{align}
Thus, finding a column of negative reduced cost for commodity $k$
amounts to computing an $s_k$--$t_k$ path that minimizes
$\sum_{e \in p} (c_e - \mu_e)$ and then subtracting the demand dual
$\pi_k$. Since all edge weights $c_e \geq 0$ and $\mu_e \leq 0$, the
modified costs $c_e - \mu_e \geq 0$, so the pricing
problem reduces to a \emph{shortest path problem in a weighted graph}.

This allows the use of efficient single-source shortest path algorithms
such as Dijkstra's method~\cite{dijkstra1959note}. In particular, when many commodities share
the same source $s$, a single run of Dijkstra's algorithm yields
shortest paths to all sinks $t_k$ with $s_k = s$, producing one pricing
candidate for each commodity simultaneously. This reduces the number of
Dijkstra runs per iteration from $|K|$ to $|S|$. The sweep itself is not new: Jones et al.\
propose it in closing \cite[p.~116]{jones1993multicommodity} and Babonneau et al.\ implement it,
one Dijkstra per source node~\cite[\S7.1]{babonneau2006active}, both running each search to
completion. It is also what makes the two formulations comparable here, since both perform the
same $|S|$ searches per iteration; \autoref{sec:implementation} adds a criterion for cutting the
sweep short.

\subsection{Tree generation}
For the tree-based formulation the procedure is:

\begin{itemize}
  \item \emph{Restricted master problem (RMP):} solve
    \eqref{eq:tree-obj}--\eqref{eq:tree} using only a subset of trees from $\cup_{s \in S} \mathcal{T}_s$.
  \item \emph{Pricing sub-problem:} for each source $s \in S$, 
  compute a shortest path tree rooted at $s$ that covers all sinks $t_k$ where $s_k = s$.
\end{itemize}

Let $\pi_s$ denote the dual variable of \eqref{eq:tree-balance} and
$\mu_e$ the dual variable of \eqref{eq:tree-cap}. The reduced cost of a
tree $\tau \in \mathcal{T}_s$ is
\begin{align}
  \bar{c}_\tau = \sum_{e \in \tau} \bar{f}^e_\tau \big(c_e - \mu_e\big) - \pi_s.
  \label{eq:tree-reduced-cost}
\end{align}
Hence, pricing reduces to finding a minimum-cost tree rooted at $s$
with respect to modified edge costs $c_e - \mu_e$. Since all modified edge
weights are non-negative, this corresponds to computing a
\emph{shortest-path tree} (SPT). A single-source shortest path
computation from source $s$ yields shortest paths $p_k$ to all sinks $t_k$ with $k \in K_s$.
Taking each $p_k$ from the predecessor structure of the same run, the union of their edges
$\tau = \bigcup_{k \in K_s} p_k$ forms a tree.
The edge flow values $\bar{f}^e_\tau$ are then computed from the paths $p_k$ by the definition in
\autoref{sec:tree}, which gives the reduced cost of the tree.

To see that the SPT minimizes the reduced cost, substitute the definition of
$\bar{f}^e_\tau$ into \eqref{eq:tree-reduced-cost} and exchange the order of summation:
\begin{align*}
\bar{c}_\tau = \sum_{k \in K_s} d_k \sum_{e \in p_k} \big(c_e - \mu_e\big) - \pi_s.
\end{align*}
Each term is minimized by choosing $p_k$ to be a shortest $s$--$t_k$ path under the modified
costs, and the SPT attains all these per-sink minima simultaneously. Since $d_k > 0$, the SPT
minimizes the demand-weighted sum and is therefore a minimum reduced cost tree.

\subsection{Implementation details}
\label{sec:implementation}

\paragraph{Bounded single-source pricing.}
We introduce the following early-termination criterion for the shared-source path pricing
described above.
For a fixed source $s$ and a set $A$ of commodities whose sinks are already settled, define
$\Pi_s := \max_{k \in K_s \setminus A} \pi_k$, the largest demand dual still outstanding. Taking
the maximum over all of $K_s$ is also valid and is the weaker form of the same test; we use the
outstanding maximum, which tightens as the search proceeds, and it is the form measured in
\autoref{sec:ablation}.
During a single run of Dijkstra on adjusted edge weights
$\tilde c_e := c_e - \mu_e$, extract-min keys are tentative distances
$d(v)$ from $s$. If the current minimum key satisfies $d(v) \ge \Pi_s$,
then for every unsettled destination $t_k$ with $k \in K_s$ we have
\[
\underbrace{d(t_k)}_{\text{shortest}} - \pi_k \;\ge\; d(v) - \pi_k
\;\ge\; \Pi_s - \pi_k \;\ge\; 0,
\]
so no path with negative reduced cost remains to be found and the pricing for source $s$
can terminate early.

The criterion carries over to tree pricing in demand-weighted form. Let $A \subseteq K_s$
denote the commodities whose sinks are settled when the minimum key is $d(v)$. Since
$\bar{c}_\tau = \sum_{k \in K_s} d_k \sum_{e \in p_k} \tilde c_e - \pi_s$ and every unsettled
sink is at distance at least $d(v)$, the run can terminate as soon as
\[
\sum_{k \in A} d_k\, d(t_k) \;+\; d(v) \sum_{k \in K_s \setminus A} d_k \;\ge\; \pi_s,
\]
because the left-hand side is a lower bound on the demand-weighted cost of any tree in
$\mathcal{T}_s$, so no tree with negative reduced cost exists.
Terminating a one-to-one shortest-path search once its minimum key exceeds the dual is standard;
what we have not found described elsewhere is the shared-source form, a single multi-target run
cut short on the maximum $\Pi_s$ over a source's demand duals, and the demand-weighted tree
analogue, which is not a rescaling of it and behaves differently (\autoref{sec:ablation}).

Both tests are stated above against zero, and in that form they are almost never satisfied: at an
optimum of the RMP every structural row carries a basic column of reduced cost zero, so a search
reaches the sinks of $s$ at a distance within rounding error of $\Pi_s$, and the test fails by
that margin. The loop, however, accepts only columns with reduced cost below a threshold $\theta < 0$, fixed
at $\theta = -10^{-3}$ below. Replacing $\Pi_s$ by
$\Pi_s + \theta$ (and $\pi_s$ by $\pi_s + \theta$ in the demand-weighted form) proves the weaker
statement the loop requires, and in that form the test is satisfied in practice: across the runs
of \autoref{sec:ablation} it terminates a median 23\% to 72\% of pricing problems early,
depending on family and formulation.

The criterion is exact in the sense required by the loop: a search it stops has already settled every
sink whose column the loop would have accepted, so a sweep returns the same columns with the
criterion on as off, which the accompanying implementation verifies by running a second, bounded
pricer against every dual vector of an unbounded run and comparing each emitted column field by
field. Under tree pricing the column is suppressed outright rather than emitted over the settled
sinks, since the master's convexity row counts trees rather than demand and would admit a tree
serving only part of its source's commodities as a complete one. The criterion does change the iteration trajectory. A truncated search never computes $\mathrm{sp}_k(c-\mu)$ for the commodities it left
unsettled, so those terms of \eqref{eq:lagrangian} fall back on the weaker value the frontier
proved and the bound is valid but looser; and the pricing filter reads the arc set recorded at a
source's last \emph{complete} price. Neither loses a column, but both can move the iteration on
which the gap exit occurs, which \autoref{sec:ablation} takes into account.

\paragraph{A* with reverse distance bounds.}
Lienkamp and Schiffer~\cite{lienkamp2024column} use A*~\cite{hart1968formal} with a precomputed
admissible heuristic to accelerate path-based pricing in the ODP setting (one source, one
destination per commodity). We extend this to the multi-target setting, where a single run from
source $s$ covers all its destinations. Let $h(v)$ be an admissible and consistent lower bound on
the distance from $v$ to any destination $t_k$ with $k \in K_s$ under the adjusted costs
$\tilde c_e = c_e - \mu_e$, and run A* with keys $f(v) = g(v) + h(v)$, $g(v)$ being the tentative
distance from $s$. The stop tests above apply with $d(v)$ replaced by $\min_v f(v)$. Exact
per-source reverse distances give the tightest $h$ but cost $O(|S||V|)$ storage, so we use a
single global bound instead: multi-source Dijkstra~\cite{dijkstra1959note} on the reverse graph
from every sink gives
\[
h(v)\;=\;\min_{k \in K}\;\text{dist}_{\text{rev}}(v,t_k),
\]
admissible and consistent for every commodity regardless of source, in $O(|V|)$ storage (see
\cite{goldberg2005computing} for landmark-based alternatives). Because $h$ is computed from the
original costs $c$ and $\tilde c_e \ge c_e$ for every $\mu \le 0$, it is admissible at every set
of duals and is computed once at start-up.

\paragraph{Partial pricing.}
A full sweep prices $|S|$ shortest-path problems, which on the larger families is itself a
substantial cost, so a sweep may stop early, either on the column quota $N$ or because the filter
described below postponed most sources. Two mechanisms ensure that no commodity is starved.
Sources are visited in round-robin order from a cursor that persists across iterations, so a sweep
cut short by the quota resumes where it stopped rather than re-pricing the same prefix; and
postponement is cleared whenever a pass returns any column. Completeness of the sweep also matters for the bound:
the bound of \eqref{eq:lagrangian} sums a term per commodity and is only valid once every source
has been visited, which is why \autoref{alg:pricing} tracks whether the sweep completed and
\autoref{alg:cgloop} consults that flag before updating the bound. Sinks left unsettled are of two
kinds: unreachable ones, whose term of the bound is dropped rather than treated as an error
(every $\mathrm{sp}_k(c-\mu) \ge 0$, so dropping it can only lower the bound), and, under the
bounded criterion only, reachable ones shown to admit no attractive column, whose term is replaced by the
frontier bound as in \autoref{alg:pricing}, so that the lower bound remains available.

\autoref{alg:search} states the per-source search and \autoref{alg:pricing} the sweep that
drives it and extracts the columns. It is the same procedure for both formulations: one
A$^\ast$ run per source, terminating once every sink of that source has been settled or the
bounded test of \autoref{sec:implementation} is satisfied, after which the predecessor structure is
walked once per commodity. Path and tree pricing
differ only in what that walk accumulates: a column per commodity, or one demand-weighted tree
column per source. The remaining quantities, $\Lambda$ and $R$, are the pricing-side terms of the
lower bound in \eqref{eq:lagrangian} below.

\providecommand{\sflag}[1]{\text{\normalfont\textsc{#1}}}  
\begin{algorithm}[htbp]
\DontPrintSemicolon
\SetKwInOut{Precompute}{Precomputed once}
\SetKwInOut{StopTest}{Stop test}
\caption{Bounded multi-target A$^\ast$ search from source $s$. One run serves every commodity
rooted at $s$; the bounded criterion of \autoref{sec:implementation} (the switch ablated in
\autoref{sec:ablation}) is tested on the frontier, so it never runs with every sink settled.}
\label{alg:search}
\Precompute{$h(v) = \min_{k \in K} \mathrm{dist}_{G^{\mathrm{rev}}}(v, t_k)$ under the
\emph{original} costs $c$, with $h(v) = \infty$ where no sink is reachable: the $O(|V|)$
dual-independent heuristic of \autoref{sec:implementation}, built once and never recomputed.}
\KwIn{source $s$; integer costs $\tilde c$ (\autoref{alg:pricing}, line~1); duals $\pi, \mu \le 0$;
flag \sflag{Bounded}}
\KwOut{settled set; distances $g$; predecessors $\mathrm{pred}$; frontier key $\hat f$; flag \sflag{Cut}}
\StopTest{With $A = \{k \in K_s : t_k \text{ settled}\}$, $\Pi_s = \max_{k \in K_s \setminus A} \pi_k$
and $D_s = \sum_{k \in K_s} d_k$, $\mathrm{Stop}(\hat f)$ is
$\hat f/\gamma \ge \Pi_s + \varepsilon_\gamma + \theta$ for path pricing and
$\bigl(\sum_{k \in A} d_k\,g(t_k) + \hat f \sum_{k \in K_s \setminus A} d_k\bigr)/\gamma
\ge \pi_s + D_s\,\varepsilon_\gamma + \theta$ for tree pricing. Keys are $\gamma$-scaled, so a
key meets a dual as $f/\gamma$, with error $\le \varepsilon_\gamma := |V|/(2\gamma)$ over a path.}
\BlankLine
$T \leftarrow \{t_k : k \in K_s\}$;\quad $r \leftarrow |T|$;\quad
$\sflag{Cut} \leftarrow \mathbf{false}$ \tcp*{did the bound stop this search?}
$g(s) \leftarrow 0$;\quad push $s$ with key $h(s)$
  \tcp*{A$^\ast$ on $\tilde c$, keys $f(v) = g(v) + h(v)$}
\While{heap $\neq \emptyset$ \KwSty{and} $r > 0$}{
  $\hat f \leftarrow$ minimum key on the heap \tcp*{the frontier, not yet settled}
  \If{\sflag{Bounded}}{
    \lIf(\tcp*[f]{no unsettled sink is reachable: heap exhaustion, not a bound hit,
      or the salvage would latch $\infty$}){$\hat f = \infty$}{\KwSty{break}}
    \lIf{$\mathrm{Stop}(\hat f)$}{$\sflag{Cut} \leftarrow \mathbf{true}$;\quad \KwSty{break}}
  }
  $v \leftarrow$ extract-min; mark $v$ settled\;
  \lIf{$v \in T$}{$r \leftarrow r - 1$ \tcp*[f]{stop once every sink of $s$ is settled}}
  \ForEach{$e = (v,w) \in \delta^+(v)$ \KwSty{with} $w$ unsettled}{
    \If{$g(v) + \tilde c_e < g(w)$}{
      $g(w) \leftarrow g(v) + \tilde c_e$;\quad $\mathrm{pred}(w) \leftarrow e$\;
      push or decrease-key $w$ with key $g(w) + h(w)$\;
    }
  }
}
\Return $(\text{settled}, g, \mathrm{pred}, \hat f, \sflag{Cut})$\;
\end{algorithm}

\begin{algorithm}[htbp]
\DontPrintSemicolon
\caption{Pricing sweep and column extraction. One search (\autoref{alg:search}) per source serves
every commodity rooted there; path and tree pricing differ only in the stop test and in lines
\ref{alg:price:pathstart}--\ref{alg:price:treeend}. $\Lambda$ and $R$ are accumulated per source
and summed in source order, so neither depends on how the batch was dispatched across threads.}
\label{alg:pricing}
\KwIn{duals $\pi, \mu \le 0$; quota $N$; cursor $\sigma$; postponed $\mathcal{Q}$; flags \sflag{Full}, \sflag{Bounded}}
\KwOut{columns $\mathcal{C}$; Lagrangian sum $\Lambda$; rounding margin $R$; flag \sflag{Swept}}
\BlankLine
\lForEach{$e \in E$}{$\tilde c_e \leftarrow \max\{0,\ \mathrm{round}((c_e - \mu_e)\,\gamma)\}$
  \tcp*[f]{$\gamma = 10^9$; integer keys}}
$\mathcal{C} \leftarrow \emptyset$;\quad $\Lambda \leftarrow 0$;\quad $R \leftarrow 0$;\quad $n \leftarrow 0$\;
\ForEach{batch $B$ of sources, taken in round-robin order from $\sigma$, or from the first
  source when \sflag{Full}}{
  $B \leftarrow$ the next sources in that order, skipping $s \in \mathcal{Q}$ unless \sflag{Full};\quad
  $n \leftarrow n + |B|$\;
  \ForEach(\tcp*[f]{the sources of a batch are priced in parallel}){$s \in B$}{
    $(\text{settled}, g, \mathrm{pred}, \hat f, \sflag{Cut}) \leftarrow$
      \autoref{alg:search} on $s$ with $\tilde c$, $\pi$, $\mu$, \sflag{Bounded}\;
    \BlankLine
    \tcp{Extraction. An unreachable sink's term is dropped ($\mathrm{sp}_k(c-\mu) \ge 0$); one
         \sflag{Cut} left unsettled is salvaged from the frontier, keeping the bound alive.}
    $\mathit{hit} \leftarrow \mathbf{false}$;\quad
    $\varphi \leftarrow \mathbf{0}$, $\bar c_\tau \leftarrow -\pi_s$ \KwSty{if} tree\;
    \lForEach{$k \in K_s$ \KwSty{with} $t_k$ unsettled \KwSty{and} \sflag{Cut}}{%
      $\Lambda \leftarrow \Lambda + d_k \max\{0,\ \hat f/\gamma - \varepsilon_\gamma\}$}
    \ForEach{$k \in K_s$ \KwSty{with} $t_k$ settled}{
      $p_k \leftarrow$ predecessor walk from $t_k$ back to $s$;\quad
      $\rho_k \leftarrow \sum_{e \in p_k} (c_e - \mu_e)$
        \tcp*{recomputed in floating point, not from $\tilde c$}
      $\Lambda \leftarrow \Lambda + d_k\,\rho_k$;\quad
      $R \leftarrow R + d_k\,|p_k|\,\epsilon_{\mathrm{LP}}$
        \tcp*{$\pi$-free; see \eqref{eq:lagrangian}}
      \uIf(\tcp*[f]{path pricing}){path}{\label{alg:price:pathstart}
        \lIf{$\rho_k - \pi_k < \theta$}{$\mathcal{C} \leftarrow \mathcal{C} \cup \{(k, p_k)\}$;\quad
          $\mathit{hit} \leftarrow \mathbf{true}$}
      }
      \Else(\tcp*[f]{tree pricing}){
        \lForEach{$e \in p_k$}{$\varphi_e \leftarrow \varphi_e + d_k$}
        $\bar c_\tau \leftarrow \bar c_\tau + d_k\,\rho_k$\;
      }
    }
    \lIf(\tcp*[f]{one tree column per source; never a partial one}){tree \KwSty{and}
      $\lnot\sflag{Cut}$ \KwSty{and} $\bar c_\tau < \theta$}{%
      $\mathcal{C} \leftarrow \mathcal{C} \cup \{(s,\varphi)\}$;\quad $\mathit{hit} \leftarrow \mathbf{true}$}
      \label{alg:price:treeend}
    \lIf(\tcp*[f]{a partial set would understate $s$'s routing}){$\lnot\sflag{Cut}$}{%
      record the walked arcs as the filter's arc set for $s$}
    \lIf{$\mathit{hit}$}{$\mathcal{Q} \leftarrow \mathcal{Q} \setminus \{s\}$}
    \lElse{$\mathcal{Q} \leftarrow \mathcal{Q} \cup \{s\}$ \tcp*[f]{postpone $s$ until its edges are cut}}
  }
  \lIf{$|\mathcal{C}| \ge N$}{\KwSty{break}
    \tcp*[f]{tested per batch, so $|\mathcal{C}|$ may overshoot $N$}}
}
$\sigma \leftarrow$ next source not yet scanned;\quad
$\sflag{Swept} \leftarrow (n = |S|)$ \tcp*{$\Lambda, R$ usable only if \sflag{Swept}}
\Return $(\mathcal{C}, \Lambda, R, \sflag{Swept})$\;
\end{algorithm}

\paragraph{Lazy addition of capacity constraints.}
In the RMP the demand constraints, \eqref{eq:path-balance} and \eqref{eq:tree-balance}
respectively, are always enforced, while the capacity constraints, \eqref{eq:path-cap} and
\eqref{eq:tree-cap}, are introduced lazily. At optimality only the binding capacity constraints
matter, since a non-binding one has a zero dual, and in optimal solutions of large
multi-commodity problems only a small fraction of edges are congested. The RMP therefore starts
with a subset (possibly empty) of the capacity rows, edge utilizations are monitored during
column generation, and the row of an edge is added when it becomes violated (row generation).
This is the \emph{active set strategy} of Babonneau et al.~\cite{babonneau2006active}, also used
by Gondzio et al.~\cite{gondzio2016primaldual}; combined with column generation it keeps the
master problem small in both dimensions.

The active set must be allowed to shrink as well as grow. An edge that is saturated early need
not stay saturated once flow has been rerouted, and a row set that only ever grows converges to
the full set of $|E|$ rows that the mechanism is intended to avoid. We therefore remove capacity rows whose dual
has been zero for five consecutive iterations. Removal is safe in the same sense that late
addition is: if the edge becomes violated again it is separated again, at the cost of one extra
iteration. The same argument applies on the column side. In master-heavy mode, where a single
iteration may add tens of thousands of columns, columns that have been out of the basis for five
iterations are purged and regenerated on demand; in pricing-heavy mode aging is disabled, because
there regeneration costs a pricing solve, which is the resource being conserved.

\paragraph{Warm start.}
Before the first solve, every pricing problem is solved once against $\pi = +\infty$ and
$\mu = 0$, and the resulting columns seed the RMP; the infinite duals make every column
attractive, so the seed contains one column per pricing problem. When slack variables sit on the
capacity rows rather than the demand rows (see below) the RMP is infeasible until every demand
row is served by at least one column, and seeding guarantees that from the first solve rather
than leaving it to a phase~1.

\paragraph{Balancing master and pricing problems.}
We propose a two-mode strategy with master-heavy and pricing-heavy operation. The bottleneck
shifts across instance families, so the mode is selected per family according to which component
dominates the runtime (\autoref{sec:cgbehaviour} reports the measurement). We do not report a
comparison against running a single mode everywhere: the design is motivated by the measured
bottleneck rather than validated against a fixed alternative, and the classification is made on
the same instances the results are reported over.

\emph{Pricing-heavy} mode is for instances on which the master RMP is the cheaper component, so
that every pricing sweep avoided is a saving. It reduces the number of sweeps:
\begin{itemize}
  \item Re-optimize the master whenever lazy cuts have been added, before attempting any pricing, so the next sweep prices against duals that already account for them
  \item Apply the pricing filter technique from \cite{lienkamp2024column}: only price sources whose previously generated columns use edges from newly added capacity constraints
  \item Clear the filter after any pricing pass that returns columns, so that postponed sources are revisited as soon as pricing succeeds again
  \item Disable column aging, since regenerating a purged column costs a pricing solve, which is the expensive resource in this mode
\end{itemize}

\emph{Master-heavy} mode is for instances on which pricing is the cheaper component. It maximizes the work
performed per master solve:
\begin{itemize}
  \item In each iteration, attempt both lazy cut separation and pricing
  \item Keep pricing until $N$ columns of negative reduced cost have been found, then stop the sweep and park the round-robin cursor where it stopped
  \item If a filtered sweep returns no column, re-price every source before terminating; this both proves optimality and, as \autoref{sec:colgen} describes below, is what makes the lower bound available
\end{itemize}

The two modes are the \texttt{pricer-heavy} and \texttt{pricer-light} presets of the accompanying
implementation, named there after the cost of the pricer rather than the resulting master
workload. The column quota $N$ is a tunable parameter. In our experiments
$N = 50\,000$ in master-heavy mode; in pricing-heavy mode the quota is not binding, as $N$ is
set to the number of pricing problems, which is the most a single sweep can return. When a sweep
returns more than $N$ columns the $N$ most negative are kept, so that a capped iteration still
passes the most attractive columns to the master rather than the first ones found. A column
is accepted only when its reduced cost is below $\theta = -10^{-3}$, chosen so that dual noise
accumulated along a path or tree cannot register as an attractive column.

\paragraph{Slack variables.}
Following standard practice in column generation~\cite{desrosiers2024branch}, we avoid
infeasibility of the master problem by adding slack variables.
When the number of demand constraints is less than the number of edges we add a slack variable to demand constraints \eqref{eq:path-balance} and \eqref{eq:tree-balance} respectively. 
Otherwise we add a slack variable to edge capacity constraints \eqref{eq:path-cap} and \eqref{eq:tree-cap}. 
The initial cost of a slack variable is the sum of all edge costs for the path-based formulation and the sum of all edge costs times the total demand for the tree-based formulation.

A single fixed cost is difficult to choose well. If it is too low, a slack stays basic at the
optimum, so the RMP objective is a penalty rather than a routing cost and no incumbent can be
recorded; if it is high enough to preclude that everywhere, it dominates the duals of the rows
it sits on, which are what pricing reads. We therefore start from the bound above and, at the end of each
iteration, multiply by ten the cost of every slack still basic with positive value, clamped at a
per-instance ceiling. Slacks that were never basic never distort the duals, a slack that is basic
only because the RMP still lacks the columns to serve its row is driven out within a few
iterations, and the clamp keeps a row that no
column can ever serve (a disconnected commodity, say) from raising its slack cost without limit
and overwhelming the objective in floating point.

\paragraph{Lagrangian lower bound.}
Relaxing the capacity constraints \eqref{eq:path-cap} with multipliers $\mu \le 0$ leaves a
problem that separates over commodities, so for any $\mu \le 0$
\begin{align}
  L(\mu) \;=\; \sum_{e \in E} u_e\, \mu_e \;+\; \sum_{k \in K} d_k\, \mathrm{sp}_k(c - \mu)
  \;\le\; z^\star,
  \label{eq:lagrangian}
\end{align}
where $\mathrm{sp}_k(c-\mu)$ is the shortest $s_k$--$t_k$ distance under the adjusted costs, which
is what the pricing sweep already computes. Taking the duals of the RMP as $\mu$ makes
\eqref{eq:lagrangian} available at no extra cost on any iteration in which every pricing problem
was solved. The structural duals $\pi$ do not appear, so the bound stays valid while a slack
variable is basic, the regime in which a bound reconstructed from the RMP objective is
meaningless because a basic slack pins $\pi$ at the slack cost. Capacity rows not yet in the
RMP carry $\mu_e = 0$ and contribute nothing, so the bound is valid
before the row set has converged. Subtracting the accumulated margin $R$ of \autoref{alg:pricing} certifies it against the
fixed-point arithmetic in which pricing runs: adjusted costs rounded to a fixed scale of $10^{9}$
and clamped at zero, so a priced path can miss the true minimum by at most the LP feasibility
tolerance $\epsilon_{\mathrm{LP}} = 10^{-4}$ per edge.
The bound is the standard Lagrangian bound of column
generation~\cite{desrosiers2024branch} applied to the capacity relaxation; we use it
both to report a gap when the time limit is reached and, in \autoref{alg:cgloop}, to stop as soon
as the incumbent is provably within $\epsilon_{\mathrm{rel}} = 10^{-4}$ of optimal.

\paragraph{Termination.}
The loop has two termination conditions, which must be distinguished when reading the results
of \autoref{sec:experiments}. The first is ordinary column generation convergence:
pricing returns nothing, separation finds no violated edge, and no slack is basic. The second is a
gap exit: the incumbent and the bound of
\eqref{eq:lagrangian} agree to within $\epsilon_{\mathrm{rel}} = 10^{-4}$ relative, at which point
further iterations are not worthwhile. Both are reported as solved, and on these instances
the two agree to well within the tolerance.

Neither exit is a proof of exact LP optimality. Pricing accepts a column only when its reduced
cost falls below $\theta$, so an empty sweep establishes that no column of reduced cost below
$\theta$ exists, and weak duality on the restricted master then bounds the gap by
$|\theta| \sum_{k \in K} d_k$ for the path formulation and, because the tree formulation's
convexity rows carry a coefficient of one rather than a demand, by $|\theta|\,|S|$ for the tree
formulation. The two exits are therefore not equally tight, and the asymmetry favours the tree formulation
by the ratio $\sum_{k \in K} d_k / |S|$. What both exits certify
in practice is the relative gap against \eqref{eq:lagrangian}, which is the quantity
\autoref{sec:experiments} reports and its acceptance check tests. An incumbent is recorded only
on iterations whose primal solution is MCF-feasible (no basic slack and no freshly violated
edge), since on any other iteration the RMP objective carries a feasibility penalty.

\paragraph{Portability across LP backends.}
The RMP is solved once per iteration by an external solver, and the study compares several,
including a GPU barrier. This has three consequences for \autoref{alg:cgloop}. All primal and dual
quantities are read immediately after the solve and before any mutation of the model, because
some backends stop returning duals once rows or columns have changed. The capacity-dual term of
\eqref{eq:lagrangian} is snapshotted against the row set the LP actually solved, which allows
the bound to be taken on an iteration that also separates. And a pure interior-point solve can
report a basic slack, or duals too central to expose an attractive column, as an artefact of a
non-vertex solution; when pricing stalls with slacks still basic we therefore request one
crossover solve, at most once per column set, and skip it on backends that offer no crossover.

\autoref{alg:cgloop} states the full loop, including where the two modes above diverge.

\providecommand{\sflag}[1]{\text{\normalfont\textsc{#1}}}  
\begin{algorithm}[htbp]
\DontPrintSemicolon
\SetKwFunction{Price}{Price}
\caption{Column generation loop: lazy capacity rows, the two operating modes, and the
Lagrangian bound that supplies the gap exit.}
\label{alg:cgloop}
\KwIn{mode $\in \{\sflag{MasterHeavy}, \sflag{PricingHeavy}\}$; iteration and wall-clock limits}
\KwOut{incumbent $\mathrm{ub}$, bound $\mathrm{lb}$, and whether the pair is certified}
\BlankLine
$N \leftarrow 50\,000$ \KwSty{if} \sflag{MasterHeavy} \KwSty{else} the most columns one sweep
  can return ($|K|$ for paths, $|S|$ for trees; either way the sweep itself runs $|S|$ searches)\;
$\mathrm{ub} \leftarrow +\infty$;\quad $\mathrm{lb} \leftarrow -\infty$;\quad
$\mathcal{Q} \leftarrow \emptyset$;\quad $\sigma \leftarrow 0$\;
add every column of \Price{$+\infty, \mathbf{0}, \infty, \sflag{Full}$}
  \tcp*{warm start, quota deliberately not applied: one column per pricing problem}
\BlankLine
\For{$t = 1, 2, \dots$}{
  \lIf{wall-clock budget exhausted}{\KwSty{break} \tcp*[f]{at the iteration boundary only}}
  solve the RMP \tcp*{exactly one LP solve per iteration}
  \lIf{not optimal}{retry once with crossover; \KwSty{break} if still not optimal}
  $z \leftarrow$ RMP objective;\quad $x \leftarrow$ primals;\quad
  $\pi \leftarrow$ structural duals;\quad $\mu \leftarrow$ capacity duals
    \tcp*{read every LP quantity before any mutation}
  $D \leftarrow \sum_{e \in E_{\mathrm{row}}} u_e\,\mu_e$
    \tcp*{snapshot against the row set just solved}
  mark each capacity row active where $|\mu_e| > \epsilon_{\mathrm{LP}}$, and each column active
    where it is basic \tcp*{ages feed the two purges below}
  \BlankLine
  \tcp{Lazy capacity rows (row generation)}
  $V \leftarrow \{e \in E : \mathrm{flow}_e(x) > u_e + \epsilon_{\mathrm{LP}}
    \text{ and } e \text{ carries no capacity row}\}$\;
  add one capacity row per $e \in V$, each with its slack column in edge-slack mode\;
  \lIf(\tcp*[f]{only then is $z$ an MCF-feasible cost}){no slack basic \KwSty{and} $V = \emptyset$}{%
    $\mathrm{ub} \leftarrow \min\{\mathrm{ub},\, z\}$}
  \lIf(\tcp*[f]{pricing filter; always on in \sflag{PricingHeavy}, optional otherwise}){filter
    \KwSty{and} $V \neq \emptyset$}{$\mathcal{Q} \leftarrow \{s \in S :
    \text{no edge of a column of } s \text{ lies in } V\}$}
  \lIf(\tcp*[f]{\sflag{PricingHeavy} digests cuts first}){\sflag{PricingHeavy}
    \KwSty{and} $V \neq \emptyset$}{\KwSty{continue}}
  \BlankLine
  \tcp{Pricing (Algorithm~\ref{alg:pricing}) on duals captured before separation}
  $(\mathcal{C}, \Lambda, R, \sflag{Swept}) \leftarrow$ \Price{$\pi, \mu, N, \lnot\sflag{Full}$}\;
  \lIf(\tcp*[f]{filtered sweep empty: re-price everything}){$\mathcal{C} = \emptyset$}{%
    $(\mathcal{C}, \Lambda, R, \sflag{Swept}) \leftarrow$ \Price{$\pi, \mu, N, \sflag{Full}$}}
  \lIf{$\mathcal{C} \neq \emptyset$}{$\mathcal{Q} \leftarrow \emptyset$
    \tcp*[f]{clear the filter, keep the cursor $\sigma$}}
  \BlankLine
  \tcp{Capacity-relaxation Lagrangian bound, valid for any $\mu \le 0$}
  \lIf{\sflag{Swept}}{$\mathrm{lb} \leftarrow \max\{\mathrm{lb},\ D + \Lambda - R\}$}
  \If{$\mathrm{ub} < \infty$ \KwSty{and}
      $0 \le \mathrm{ub} - \mathrm{lb} < \epsilon_{\mathrm{rel}} \max\{1, |\mathrm{ub}|\}$}{
    \Return $\mathrm{ub}$, certified \tcp*{gap exit; $\mathcal{C}$ is discarded unadded}
  }
  \lIf{$|\mathcal{C}| > N$}{keep the $N$ most negative reduced costs}
  \BlankLine
  \If{$\mathcal{C} = \emptyset$}{
    \lIf(\tcp*[f]{no attractive column, no violation, no slack}){$V = \emptyset$
      \KwSty{and} no slack basic}{\Return $\mathrm{ub}$, certified}
    request one crossover solve next iteration (at most once per column set)\;
    multiply the cost of every basic slack by $10$, clamped at its ceiling\;
    $\mathcal{Q} \leftarrow \emptyset$;\quad $\sigma \leftarrow 0$;\quad \KwSty{continue}\;
  }
  multiply the cost of every basic slack by $10$, clamped at its ceiling\;
  purge columns idle for $\ge 5$ iterations \tcp*{\sflag{MasterHeavy} only; aging is off otherwise}
  purge capacity rows non-binding for $\ge 5$ iterations\;
  add $\mathcal{C}$ to the RMP\;
}
\Return $\mathrm{ub}$ if finite else $\mathrm{lb}$, not certified\;
\end{algorithm}

\paragraph{Contribution of the individual techniques.}
The contribution of each technique in isolation is measured for the bounded stopping criterion
only, by the paired comparison in \autoref{sec:ablation}; no leave-one-out measurement is
reported for the others. The criterion is ours, so no published measurement covers it, and it
can be removed without side effects: switching it off changes when a search stops, not what the loop
does with the result, so both arms of the comparison finish and can be compared directly.
The other techniques are adopted, and their leave-one-out measurements have been published by the
authors from whom we adopt them: the lazy capacity mechanism is the active-set strategy of Babonneau et
al.~\cite{babonneau2006active}, and Lienkamp and Schiffer report their pricing filter cutting
computation time by up to 60\% and their A*-based pricing a further 90\% on top of
it~\cite{lienkamp2024column}, measured on the same intermodal family we use here; our
contribution to that pair is the extension of the A* bound to the multi-target tree setting,
not the techniques themselves. Our own measurements also bound what any such ablation could
show: on the grid, planar and transportation families pricing is at most 5.3\% of runtime
(\autoref{sec:cgbehaviour}), so no pricing-side technique can change their wall-clock time by more than that, whatever its
effect on pricing itself. The intermodal family is the exception, with the balance reversed, and
the two-mode strategy therefore selects pricing-heavy operation for it.

\section{Experiments}
\label{sec:experiments}

We evaluate the path-based and tree-based formulations on 44 instances from four families, using
the source-based formulation solved directly as a baseline.

The comparison of two decomposition schemes depends on the linear programming solver used for
the restricted master problems: each barrier code has its own scaling behaviour, numerical
tolerances and failure modes.
We therefore run the entire study under five LP configurations rather than one, spanning
open-source and commercial codes and both CPU and GPU execution:

\begin{itemize}
  \item \textbf{HiGHS 1.15.1}~\cite{huangfu2018parallelizing} with the HiPO interior-point
        solver~\cite{zanetti2025factorisation} (open source, CPU),
        carrying a local patch to HiPO's IPX refinement step: without it HiPO discards a
        near-optimal solution as an internal error on the ill-conditioned path masters,
        which then fail rather than solve.
  \item \textbf{MOSEK 11.0.30}~\cite{mosek} interior point (commercial, CPU).
  \item \textbf{cuOpt}~\cite{cuopt} barrier (open source, GPU). The column generation runs use a fork
        (reporting version 26.8.0) that adds an incremental C API for model modification.
        Stock cuOpt exposes no such interface, so the restricted master would be rebuilt
        from scratch at every iteration, a serious performance degradation and not
        representative of what the solver can do on this workload. The compact baseline,
        a single solve that never mutates its model, uses the stock 26.2.0 release.
  \item \textbf{COPT 8.0.1}~\cite{ge2022copt} barrier, in CPU mode and in GPU mode (commercial); the two modes
        differ only in barrier execution, which isolates the effect of the accelerator from
        the effect of the solver.
\end{itemize}

Two of the five are therefore not off-the-shelf builds, and neither deviation is visible in
the version the library reports. Both are open source and pinned: the HiGHS patch is applied
at configure time from the implementation repository, and the cuOpt fork is
\url{https://github.com/spoorendonk/cuopt} at commit
\texttt{8ea7a033addb5594532f8f60b33b99067459c853}, cited by commit because the
\texttt{delta-api} branch is a moving reference. Each configuration solves both the restricted
master problems of the two column generation formulations and the compact source-based model of
the baseline, so the tree-versus-path comparison is like-for-like within a single LP code
throughout. The one exception is the compact baseline under cuOpt, which runs on the stock
26.2.0 release rather than the fork; the fork also carries a barrier fix for a
ranged-constraint miscomputation with a large finite bound, so the cuOpt source-based figures
should be read with that caveat. Crossover is disabled throughout (apart from the
stall-recovery solve \autoref{sec:colgen} requests at most once per column set), as is presolve
on the restricted master problems; on the compact baseline presolve is left at each solver's
default, since a direct monolithic solve benefits from it and forcing it off would
handicap the baseline.

All experiments were performed on a machine with an AMD Ryzen 9 3950X 16-Core Processor,
128\,GB of RAM and an NVIDIA GeForce RTX 5070 Ti, running Ubuntu 24.04. The algorithms were
implemented in C++, using up to 32 threads.

Throughout the experiments we set a relative optimality tolerance of $10^{-4}$ and a timeout of
7200 seconds. The grid, planar and transportation families run in master-heavy mode and the intermodal
family in pricing-heavy mode (\autoref{sec:colgen}); the bounded stopping criterion follows
the same split, enabled on the intermodal family and
disabled on the other three. This choice is based on the measurements of \autoref{sec:ablation}: the criterion reduces
pricing time on every family, but only on intermodal is pricing a large enough share of the
runtime for the reduction to affect the total.

A run counts as solved only if the LP backend reports an optimal status \emph{and} the objective
matches the reference value to $10^{-3}$ relative. That threshold is looser than the $10^{-4}$
the loop targets, but 419 of the 420 solved column generation runs land within $10^{-4}$, and
the exception is a path-based run at $3.6 \times 10^{-4}$, so the looser threshold does not favour the
tree-based formulation. A run whose objective matches but which terminates at the time limit
counts as unsolved (three runs of the compact model do), and runs that fail either
test enter the aggregate statistics at the time limit. Because the limit is enforced at
iteration boundaries, a recorded wall-clock time can exceed it slightly, and we report the limit
rather than the overrun.

The reference values are external where one exists: the grid and planar families are published
with the instances, and the intermodal optima come from the generator. For the transportation
family no external reference applies, because the demand coefficients differ from those of
\cite{babonneau2006active}; there the reference is our own, the value the five backends agree
on, and on Austin and Philadelphia, which the path-based formulation reaches under no backend,
the agreement is across the five tree-based runs alone. The compact baseline is solved at the
same $10^{-4}$ barrier tolerance with crossover disabled as the restricted master problems, so
the advantage reported over it in \autoref{sec:performance} is not an artefact of holding the
baseline to a stricter standard.

\subsection{Instances}
\label{sec:instances}
We consider 44 instances across four families: 15 grid, 10 planar, 9 transportation networks and
10 intermodal transport networks. The grid and planar instances are from the MCF instances
benchmark \cite{mcf-instances}, the transportation networks from the transportation networks
benchmark \cite{transportationnetworks2025}, and the intermodal transport networks from the
benchmark of \cite{lienkamp2024column}.

The transportation networks are non-linear maximum flow MCF problems, and as in
\cite{babonneau2006active} we convert them to min-cost MCF problems: ``free flow time'' is the
unit edge cost, and the demands are divided by a coefficient, increased until the problem
becomes feasible. \autoref{tbl:coefs} in \autoref{sec:appendix} lists the coefficients.
Babonneau et al.~\cite{babonneau2006active} report the transportation set at a relative gap of
$10^{-5}$ and a separate low-precision run at $10^{-3}$, but we were not able to reproduce their
result with their coefficient values on the current version of the instances
\cite{transportationnetworks2025}, so our coefficients and optimal objective values differ. On
ChicagoSketch, ChicagoRegional and Philadelphia the difference is small (2.4, 4.1 and 7.0
against their 2.5, 6 and 7), but on Winnipeg and Barcelona it is three orders of magnitude
(2000.0 and 5050.0 against 2.7 and 3). The coefficient divides the demands, so a larger value
leaves a less congested and therefore easier linear program: our transportation instances are
not the same problems theirs were, and the comparison in \autoref{sec:performance} is
indicative rather than controlled.

The intermodal transport networks are derived from the Munich public transport network,
generated with the code of \cite{lienkamp2024column} with seed 0; preliminary tests showed seeds
1 to 9 to be computationally equivalent, so we consider seed 0 only, and we filter out trips
that are not connected in the generated network. The generator also produces a SUBWAY family,
which we exclude because no published reference optimum is available for it, leaving the five
BUS and five SBT instances.

\autoref{tbl:familytraits} shows the characteristics of the instance families, and
\autoref{tbl:instancedetails} in \autoref{sec:appendix} the full per-instance figures. Grid and
planar are relatively small but most of their vertices are commodity sources; the transportation
networks have a large number of commodities with many shared sources; the intermodal networks
are large in vertices and edges but have relatively few commodities and no shared sources.

The intermodal family requires justification, since source aggregation gains nothing there:
every commodity has its own source, so a shortest-path tree \emph{is} a path and the two
formulations reduce to the same model. It is retained as the control for the study, the family
on which the two formulations are known a priori to be identical, so
any difference the measurements show there is measurement error rather than a
formulation effect; \autoref{sec:cgbehaviour} and \autoref{sec:planar2500} use it in this role.
Without such a family the differences reported on the other three could not be separated from
artefacts of the implementation.

\begin{table}[htbp]
\centering
\begin{adjustbox}{max width=\textwidth}
\begin{tabular}{lrrrrrrr}
\toprule
Family & Nodes & Edges & Commodities & Sources & Variables & Constraints & Non-zeros \\
\midrule
Grid & 523 & 2,011 & 4,967 & 461 & 1.6e+06 & 4.1e+05 & 4.7e+06 \\
Planar & 551 & 2,932 & 12,555 & 551 & 4.3e+06 & 8.3e+05 & 1.3e+07 \\
Transportation & 10,833 & 27,109 & 944,610 & 1,113 & 5.0e+07 & 2.0e+07 & 1.5e+08 \\
Intermodal & 177,220 & 1,894,289 & 22,195 & 22,195 & 6.4e+10 & 4.6e+09 & 1.9e+11 \\
\bottomrule
\end{tabular}
\end{adjustbox}
\vspace{1em}
\caption{Instance family characteristics (mean values). Number of variables, constraints, and non-zeros are based on the source-based formulation.}
\label{tbl:familytraits}
\end{table}

\subsection{Performance}
\label{sec:performance}
\begin{figure}[htbp]
  \centering
  \includegraphics[width=\textwidth]{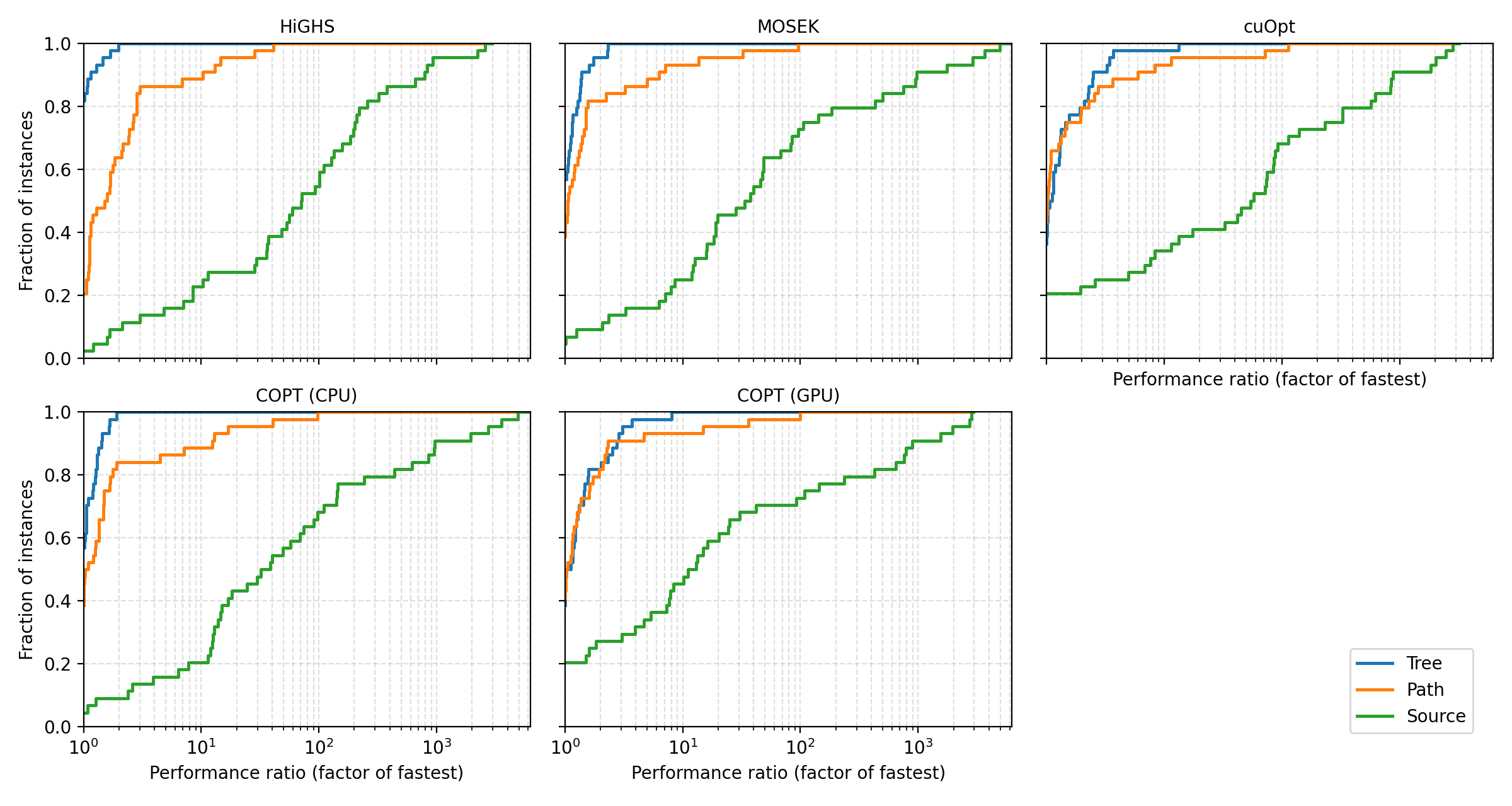}
  \caption{Performance profile over all 44 instances, one panel per LP backend: fraction of instances solved (vertical axis) within a given factor of the fastest approach's runtime (horizontal axis)}
  \label{fig:performance-profile}
\end{figure}

\autoref{fig:performance-profile} shows the performance profile of the three approaches, with one
panel per LP backend. The aggregate ordering is the same in every panel: the tree-based
formulation's curve dominates, the path-based formulation's lies below it, and the direct solve of
the source-based model is last, matching the shifted geometric means of \autoref{tbl:sqmstats}. The reproduction of this ranking by four
independent barrier codes in five configurations (open source and commercial, CPU and GPU) is
the central robustness check of this study: the ordering is not an artefact of any single interior point
implementation. All five configurations are interior point methods with crossover disabled. No simplex master is
among them, by decision rather than omission: preliminary tests found the simplex far slower
than interior point on masters of this size, and because the loop adds capacity rows and columns
in the same iteration, the warm start that is the simplex's advantage in column generation is
largely lost.
Neither of the two prior large-scale studies uses a simplex master either: Babonneau et
al.~\cite{babonneau2006active} relax the capacity constraints and solve the Lagrangian dual
with Proximal-ACCPM, an analytic-centre cutting-plane method, and Gondzio et
al.~\cite{gondzio2016primaldual} solve the restricted master with a primal-dual interior point
method, citing the dual oscillation and degeneracy of simplex-based masters. Babonneau et al.\
also record the earlier consensus that a fast simplex code was the fastest choice for small and
medium instances, which is the regime Jones et al.\ worked in.

The profiles also show that the ordering holds in aggregate and not instance by instance. Under
cuOpt and COPT in GPU mode the path-based formulation is faster on a majority of the
individual instances the two both solve (24 of 39 and 27 of 42 respectively), and it loses on the
aggregate only because the margins by which it wins are small while those by which it loses are
large: where the tree-based formulation is faster, it is faster by up to two orders of magnitude.
The aggregate ordering, not the instance-level one, is invariant across backends.

The tree-based formulation solves 43 or 44 of the 44 instances depending on the backend; its only
failure is planar2500, and only under HiGHS and cuOpt. The path-based formulation solves between
38 and 42, failing on 17 of the 45 transportation runs and on planar2500 under HiGHS. Those
failures are concentrated: on Austin and Philadelphia the path-based formulation reaches the
required accuracy under \emph{no} backend; on Sydney it fails under four backends of five, on
ChicagoRegional under two, and on Birmingham and planar2500 under one. The direct solve of the
source-based model solves between 24 and 27.

Nine of the 44 instances have no source-based entry. The baseline runs in three steps:
build the compact model in memory, stream it to a gzipped MPS file, and hand that file to the
vendor's command-line barrier. The source-based model has $|S||E|$ variables with at most three
non-zeros each, so on the intermodal family, where every commodity has its own source and
$|S| = |K|$, both counts are correspondingly large. On nine of the ten intermodal instances the first step is
not completed, and cannot be completed on this machine: at the twelve bytes per non-zero of any
sparse representation their constraint matrices alone would occupy from 164\,GB (SBT-6255-0) to
9.4\,TB (SBT-56295-0) against the 128\,GB of RAM, before any solver is invoked. The tenth,
BUS-2632-0, and Sydney were exported and are solved by no backend: on eight of their ten runs
the kernel terminates the process with its resident set at 105\,GB, and on the other two COPT's
MPS reader aborts at 60\,GB before a model reaches the barrier. We treat the nine unbuilt
instances as unsolved by the source-based approach (failing to build the model is a failure to
solve it), which is in itself an argument for decomposition. Detailed results for the
instances are shown in \autoref{tbl:runtime} in \autoref{sec:appendix}.

\begin{figure}[htbp]
  \centering
  \includegraphics[width=\textwidth]{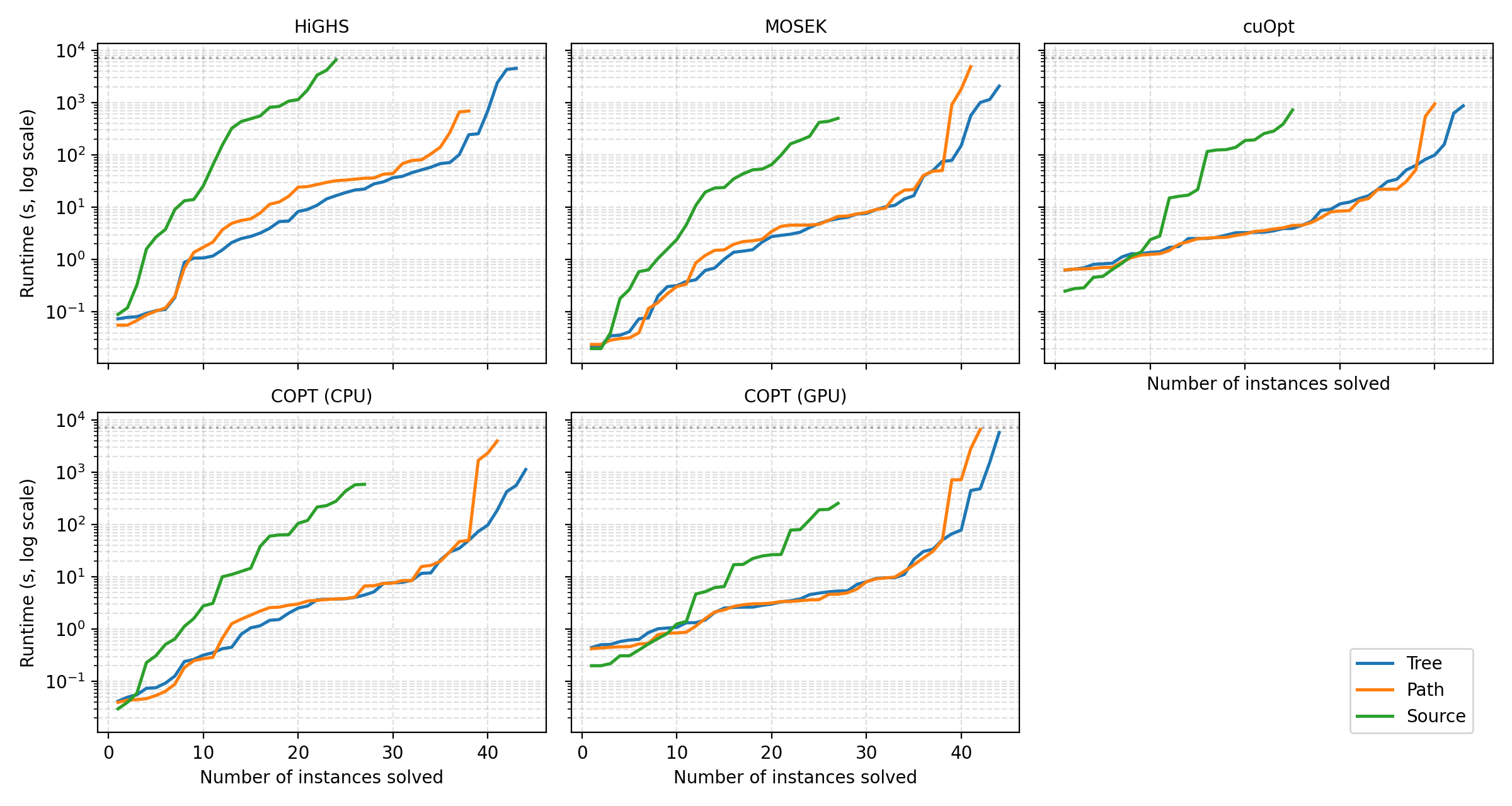}
  \caption{Cactus plot of runtime, one panel per LP backend: number of instances solved (horizontal axis) within a given runtime (vertical axis, log scale)}
  \label{fig:cactus-plot-runtime}
\end{figure}

\autoref{fig:cactus-plot-runtime} shows the corresponding cactus plot. In every panel the
tree-based curve extends furthest to the right and runs below the path-based curve in the tail,
where the large transportation instances sit.

\autoref{tbl:sqmstats} reports the shifted geometric mean runtime per backend. Means are shifted
by 10 seconds to reduce the impact of very small runtimes, following standard practice in solver
benchmarking \cite{achterberg2007constraint}. The tree-based formulation is faster than the
path-based one under every backend, by a factor between 1.42 and 1.93, and faster than the direct
solve by a factor between 15 and 28. The tree-versus-path factor is conservative: 18 of the 20 failed column generation runs are
path-based and enter the mean at the 7200\,s limit rather than at the time they would have
needed, so the true gap is wider than 1.42 to 1.93. The factor over the compact baseline is, conversely, inflated by non-completion: of the 220
compact cells, 90 enter at the limit, and only eight of those are timeouts, the rest being 45 unbuilt models and 37
crashes or early stops. Restricted to the instances the compact model actually solves under a
given backend, the factor falls to between 5 and 25. Across backends, the fully open-source
configuration (HiGHS with HiPO) is roughly 2.8 times slower on the mean than
the fastest commercial configuration, so a reader reproducing this work without commercial
licences should expect the same ordering at a uniformly higher cost.

The `Combined' row of \autoref{tbl:sqmstats} bounds the gain available from choosing the
formulation per instance. Taking
the better of the two runs on each instance solves all 44 instances under four of the five
backends, because planar2500, the only instance the tree-based formulation ever misses, is
solved by the path-based formulation under cuOpt; HiGHS leaves it unsolved either way. On
runtime, picking the better formulation for every instance saves between 2\% and 19\% of the
tree-based mean (2\% under HiGHS and under COPT in CPU mode, 6\% under MOSEK, 13\% under COPT in
GPU mode and 19\% under cuOpt), whereas the path-based formulation alone is between 1.6 and 2.0
times that same bound. The row is a bound rather than a method, since choosing the formulation
per instance presupposes knowing the outcome, and it shows that the gain over always choosing
the tree-based formulation is small.

The largest single-instance gap between the two formulations, among instances that both of them
solve, is Sydney under COPT in GPU mode: 66 seconds for the tree-based formulation against 6563
seconds for the path-based one, a factor of 99. Restricting to such instances understates the
difference, since the cases where the path-based formulation is furthest behind are those
where it does not finish.

\begin{table}[htbp]
\centering
\begin{tabular}{llrrr}
\toprule
LP backend & Formulation & SGM time (s) & Scaled & Solved \\
\midrule
HiGHS & Tree & 31.8 & 1.00 & 43 of 44 \\
 & Path & 57.9 & 1.82 & 38 of 44 \\
 & Combined & 31.2 & 0.98 & 43 of 44 \\
 & Source & 888.0 & 27.91 & 24 of 44 \\
\addlinespace
MOSEK & Tree & 14.0 & 1.00 & 44 of 44 \\
 & Path & 22.4 & 1.59 & 41 of 44 \\
 & Combined & 13.2 & 0.94 & 44 of 44 \\
 & Source & 285.2 & 20.31 & 27 of 44 \\
\addlinespace
cuOpt & Tree & 15.2 & 1.00 & 43 of 44 \\
 & Path & 21.9 & 1.44 & 40 of 44 \\
 & Combined & 12.3 & 0.81 & 44 of 44 \\
 & Source & 380.2 & 24.94 & 25 of 44 \\
\addlinespace
COPT (CPU) & Tree & 11.5 & 1.00 & 44 of 44 \\
 & Path & 22.1 & 1.93 & 41 of 44 \\
 & Combined & 11.2 & 0.98 & 44 of 44 \\
 & Source & 293.0 & 25.58 & 27 of 44 \\
\addlinespace
COPT (GPU) & Tree & 14.1 & 1.00 & 44 of 44 \\
 & Path & 20.0 & 1.42 & 42 of 44 \\
 & Combined & 12.3 & 0.87 & 44 of 44 \\
 & Source & 217.7 & 15.45 & 27 of 44 \\
\bottomrule
\end{tabular}
\vspace{1em}
\caption{Shifted geometric mean runtime (seconds, shift 10\,s) over the 44 instances, per LP backend. Scaled columns are relative to the tree-based formulation for that backend. Runs that did not reach the required accuracy enter the mean at the time limit (7200\,s), as do the nine instances for which no compact source-based model was ever built, and so none solved. The `Combined' row takes, for each instance, the better of the tree- and path-based runs; it is not a method that can be run without knowing the outcome in advance, but an upper bound on the gain from choosing the formulation per instance.}
\label{tbl:sqmstats}
\end{table}

\paragraph{Comparison to other work.}
We compare our implementation against several recent approaches from the literature. These
comparisons reflect our full system (formulation choice plus the engineering techniques of
\autoref{sec:colgen}) rather than the formulation in isolation, and the hardware differs from
the cited work, so they are indicative rather than controlled.

For planar2500, the largest instance solved in \cite{gondzio2016primaldual}, our path-based
implementation takes between 544 seconds (cuOpt) and 1683 seconds (COPT in CPU mode) against
the 7102 seconds reported there (their optimality tolerance is tighter, $10^{-5}$ against our
$10^{-4}$, and their hardware older), between 4.2 and 13 times faster under the four backends
that solve it; under HiGHS it terminates after 56 seconds without reaching the required
accuracy. This comparison concerns the implementation rather than the formulation; planar2500 is
discussed as a formulation question in \autoref{sec:planar2500}.

The transportation network instances from \cite{babonneau2006active} include cases that were solved only
with their active-set strategy: without it, neither Proximal-ACCPM nor ACCPM solved
ChicagoRegional or Philadelphia, and with it they report 14684 and 3125 seconds
respectively, at a relative gap of $10^{-5}$ on a 2.8\,GHz Pentium~IV. The tree-based
formulation solves all nine instances to $10^{-4}$ accuracy under every backend, no
transportation instance taking more than 559 seconds under COPT in CPU mode; the most expensive instance
across backends is Austin, from 482 seconds under COPT in GPU mode to 4518 under HiGHS, and the
larger Sydney network is solved in 63 to 255 seconds. This is not a direct comparison, since
the demand coefficients differ from theirs (\autoref{sec:instances}).

For the intermodal transport networks, we compare against \cite{lienkamp2024column} who solved
their largest instance (SBT-56295) in 2856 seconds on an Intel Core i9-9900 (3.1 GHz). That
figure is a mean over ten randomly generated instances of that size, and it is the time for their
full price-and-branch procedure including the integer step, where we report the linear relaxation
alone.
Accounting for processor speed differences based on PassMark multi-core CPU benchmark scores,
their solution time normalizes to approximately 1200 seconds on our AMD Ryzen 9 3950X (3.5 GHz).
Our implementation solves the same instance in 50 to 81 seconds depending on the backend and
formulation (50 to 72 under the tree-based formulation), more
than an order of magnitude faster.

\subsection{Memory}
\label{sec:memory}

\begin{figure}[htbp]
  \centering
  \includegraphics[width=\textwidth]{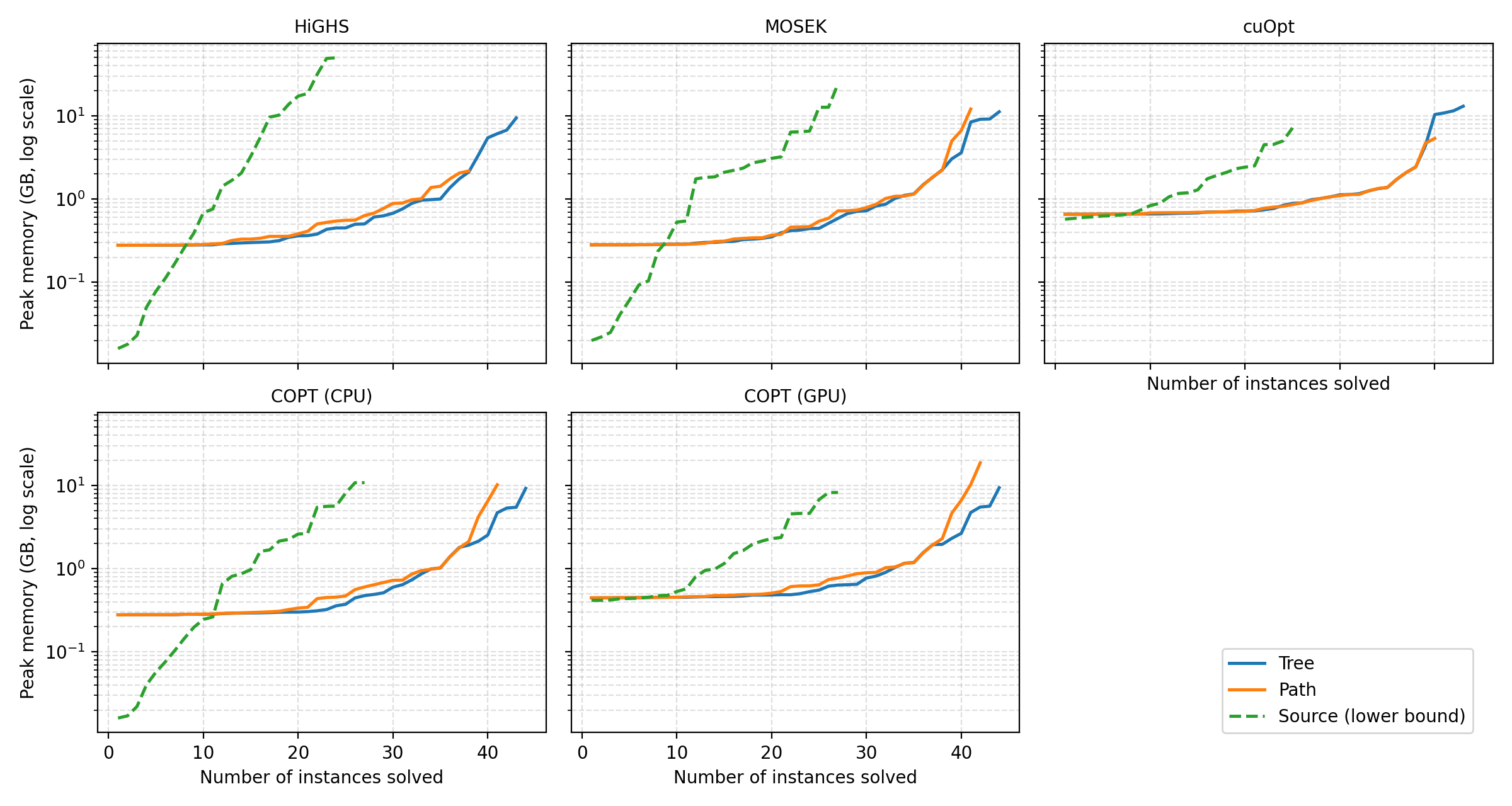}
  \caption{Cactus plot of peak resident memory, one panel per LP backend: number of instances
  solved (horizontal axis) within a given peak memory (vertical axis, log scale). Only solved
  runs are plotted, so each curve ends at the number of instances that approach solved under
  that backend. The source-based curve is a lower bound, dashed for that reason: its figure is
  the peak of a barrier run stopped after three iterations (model read, presolve,
  factorization and three iterations), measured in a separate execution from the timed solve,
  so it understates what the full solve needs}
  \label{fig:cactus-plot-memory}
\end{figure}

Since the tree-based master problem is the smaller of the two, it is expected to require less
memory. It does, but the difference is much smaller than the runtime difference. Over the 201 (instance,
LP backend) cells that both formulations solved, the geometric mean ratio of path to tree peak
resident memory is 1.09, with a median of 1.00, a minimum of 0.71 and a maximum of 2.86. By family, it is
1.00 on intermodal, where the two formulations are the same model, 1.04 on grid, 1.14 on planar
and 1.37 on transportation, the family with the largest runtime advantage, and it reaches 2.0 on
Birmingham and Sydney. \autoref{fig:cactus-plot-memory} shows the same per
backend: the two curves are indistinguishable over the bulk of the suite and separate only in
the tail.

Memory is the limiting factor in the runs the path-based formulation does not finish. Of the 20 runs that
failed to reach the required accuracy, 18 are path-based, and the largest memory readings in the
study are among them: Sydney under the path-based formulation and MOSEK is the only one of
the 440 runs that did not exit cleanly, terminated by the kernel at 95.8\,GB of resident memory
on the 128\,GB machine, where the tree-based formulation solved the same instance under the same
backend in 76 seconds with a peak of 11.2\,GB. No run of either formulation that solved its
instance exceeded 18.7\,GB, and no tree-based run anywhere exceeded 13.1\,GB. Unfinished runs
are reported at their measured peak rather than at the machine limit, and \autoref{tbl:memory}
in \autoref{sec:appendix} marks them.

Peak resident memory is a whole-process figure; on this suite it is dominated by the master.
Everything on the pricing side is sized by the network: one copy of the graph in compressed
sparse row form shared by every pricing thread, the reverse distance bounds of
\autoref{sec:colgen}, one value per vertex, and one Dijkstra workspace of $O(|V|)$ per thread,
with nothing replicated per commodity. The master is sized by the row and column counts of
\autoref{tbl:cgbehaviour}, and on the instances where memory is limiting the two differ by orders of
magnitude: Sydney has 33\,113 vertices and 75\,379 arcs, against a path-based master
carrying a demand row and a seed column for each of its 3.3 million commodities. A 95.8\,GB
peak cannot be attributed to structures sized by the former, so the ratios
above compare the two formulations' masters.

The compact model's memory is known only as a lower bound. The baseline runs through each
vendor's command-line solver, and the peak we have for every cell is that of a barrier run
stopped after three iterations in a separate execution, which understates the full solve. Even
so, on the 130 cells the compact model solves that bound is a geometric mean 2.5 times the
tree-based formulation's full-run peak, above it on 90 of the 130 and 3.5 times on the
transportation family; and on the nine instances it never builds, the constraint matrix alone
is 164\,GB to 9.4\,TB (\autoref{sec:performance}). Both decompositions solve instances for
which the undecomposed model does not fit in memory.

\subsection{Comparison of the path- and tree-based formulations}
\begin{figure}[htbp]
  \centering
  \includegraphics[width=0.8\textwidth]{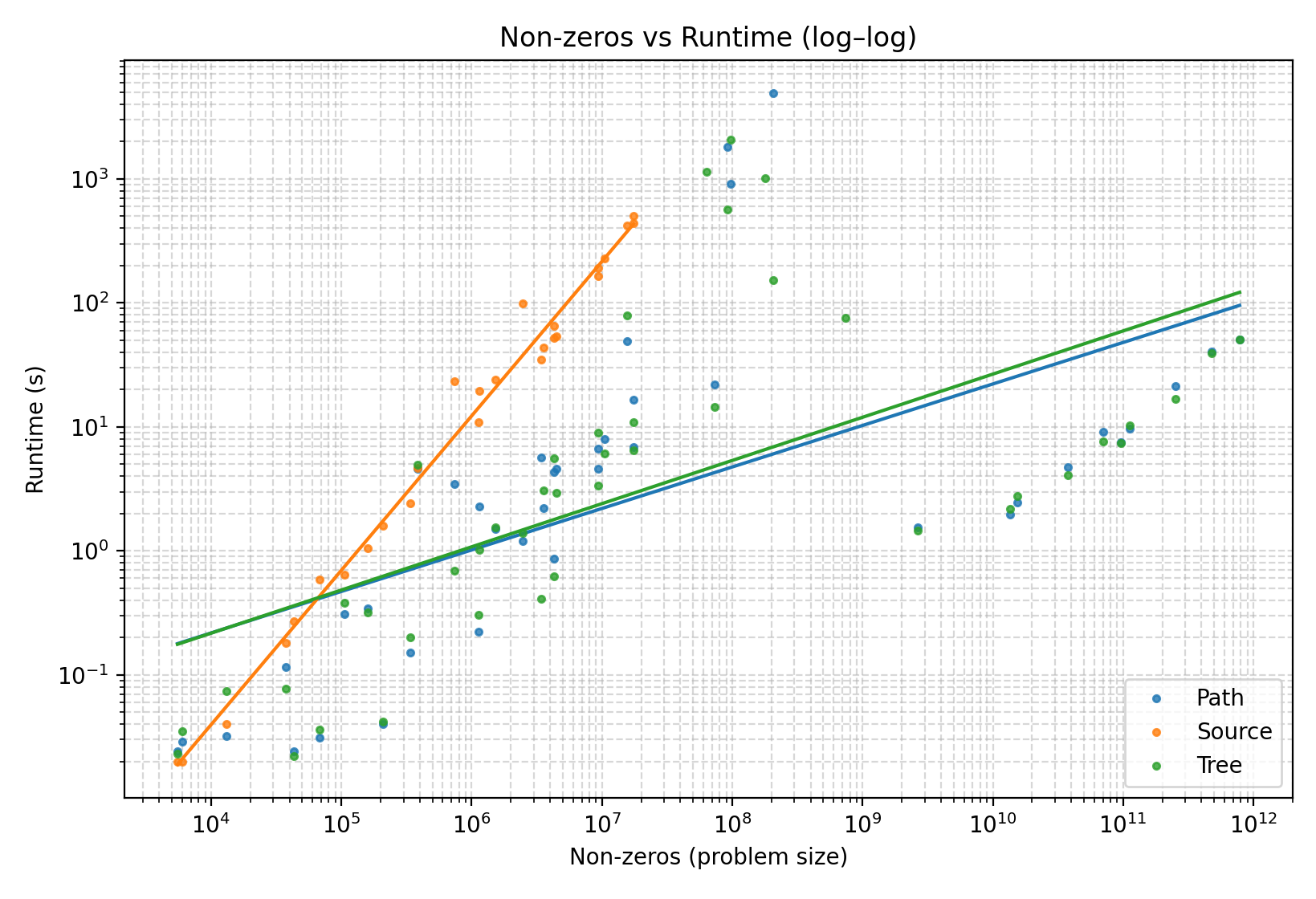}
  \caption{Runtime vs size in number of non-zeros of the source-based formulation, under the
  MOSEK backend. Runs that reached the time limit are omitted rather than plotted at it, which on
  this backend drops three path-based points and no tree-based ones.}
  \label{fig:size-vs-time}
\end{figure}

\autoref{fig:size-vs-time} plots runtime against instance size, measured as the number of
non-zeros of the source-based LP, under the MOSEK backend. The compact solve's runtime grows far
more steeply with size than either decomposition's, and its points stop at $10^{7}$ non-zeros
where the decompositions continue past $10^{11}$. The path- and tree-based fits are
near-parallel: the two formulations scale alike in problem size, and they differ by a
factor determined by source sharing rather than by growth rate. The intermodal instances, at the
far right, are those on which the two coincide (\autoref{sec:instances}).

\begin{figure}[htbp]
  \centering
  \includegraphics[width=0.8\textwidth]{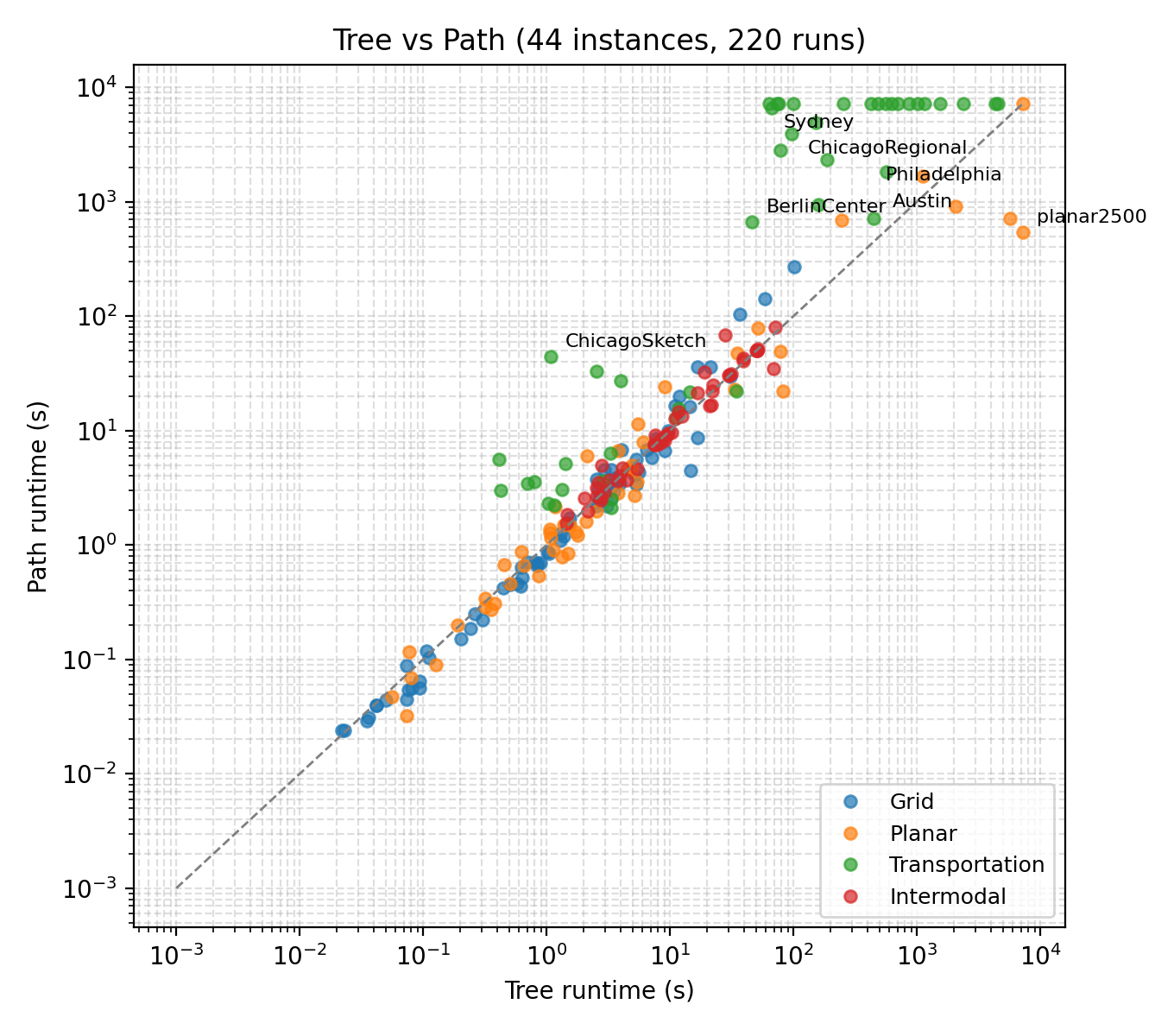}
  \caption{Instance-wise comparison of tree-based (horizontal axis) and path-based (vertical axis) runtimes on log scales, one point per instance and LP backend; points above the diagonal are solved faster by the tree-based formulation}
  \label{fig:scatter-plot-path-tree}
\end{figure}

\autoref{fig:scatter-plot-path-tree} compares the two formulations instance by instance, with one
point per instance and LP backend. Points cluster by instance, but 28 of the 41 instances that
both formulations solve under at least two backends have a different winner under at least one
of them. The flips are concentrated where the two formulations are close (sub-second grid and
planar instances, and the intermodal family), and the sign is stable where the effect is
large: on the transportation family the tree-based formulation is faster in 25 of the 28 cells
both solve, the exceptions being BerlinCenter under cuOpt and Winnipeg under cuOpt and COPT in
GPU mode, all three among the family's smallest instances. The transportation instances, with
their many shared sources, lie well above the diagonal, and planar2500 is the clearest case below it.

\begin{figure}[htbp]
  \centering
  \includegraphics[width=0.8\textwidth]{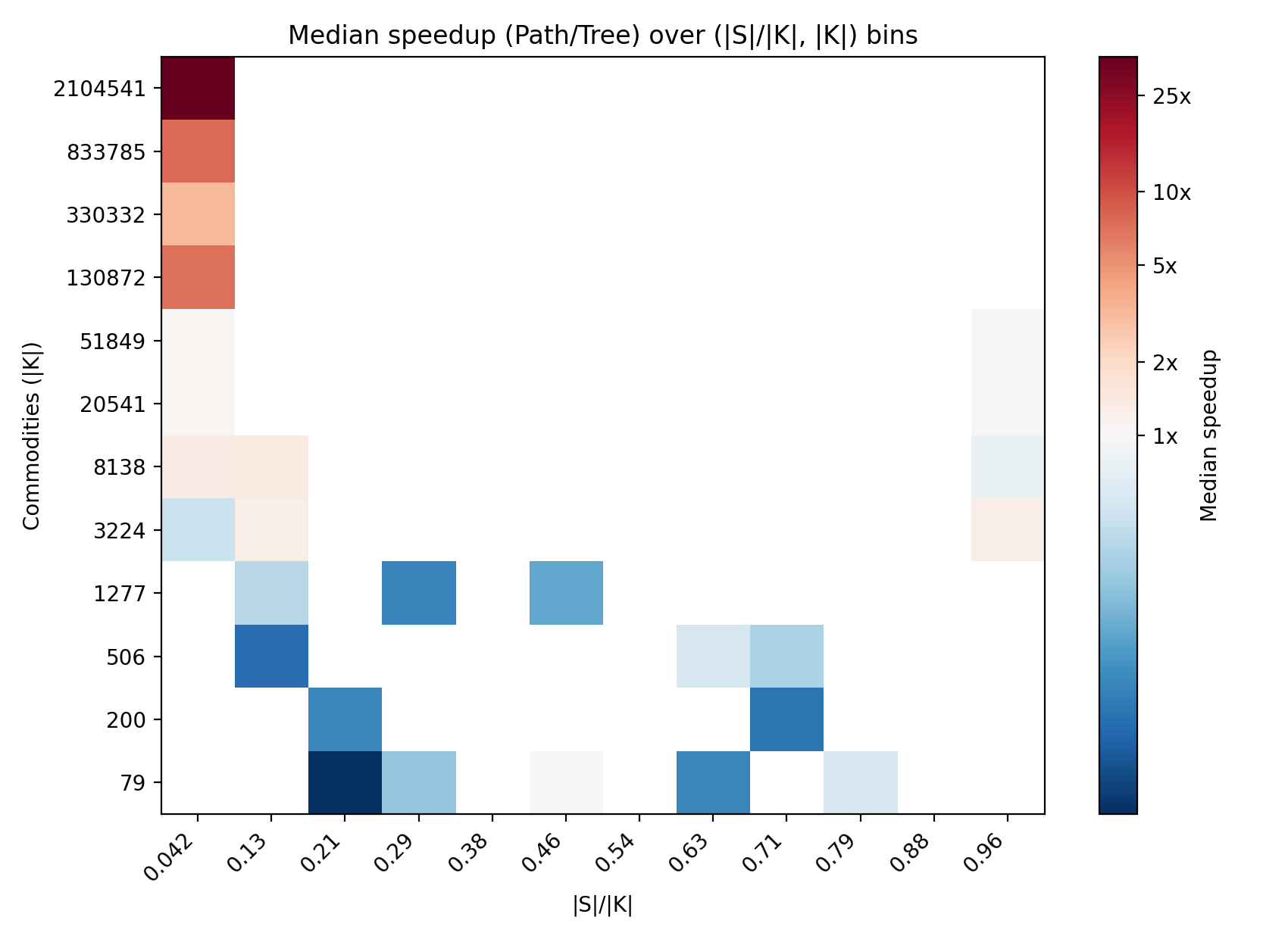}
  \caption{Median speed-up of the tree-based over the path-based formulation (ratio of path to tree runtime) as a function of the source ratio $|S|/|K|$ (horizontal axis) and the number of commodities $|K|$ (vertical axis)}
  \label{fig:heatmap-path-tree}
\end{figure}

The heatmap in \autoref{fig:heatmap-path-tree} characterizes empirically when to prefer each
formulation. It plots the median speed-up of tree over path, pooled over all five backends, as a
function of $|K|$ and the source ratio $|S|/|K|$: the tree-based formulation wins for large
$|K|$ and small $|S|/|K|$ (red) and loses in the opposite corner (blue), consistent with the
$O(|S|)$ against $O(|K|)$ demand constraint count. The instances available to us do not,
however, separate $|K|$ from $|S|/|K|$. The large-$|K|$ instances in the
public benchmarks are all heavily source-shared, and the only instances with $|S|/|K| \approx 1$
at scale are the intermodal ones, where the two formulations are the same model by construction,
so the heatmap describes the regimes these benchmarks populate rather than an identified boundary
between them. Within those regimes it supports a practical rule: prefer the source-based tree
decomposition when $|S| \ll |K|$ and $|K|$ is large. The source ratio alone is not sufficient.
planar2500 has the \emph{lowest} $|S|/|K|$ in the planar family and is the instance on which
the tree-based formulation loses, under three of the five backends, by up to a factor of
eight among the cells both formulations solve and by 13.2 against cuOpt, where the tree-based run
does not finish; \autoref{sec:planar2500} analyses it and identifies a second condition for the
rule, that a shortest-path tree from a source be sparse relative to the master's row count,
which fails when $|S|$ approaches $|V|$. Outside the regime where both conditions hold the two
formulations are within a small factor and the winner is not stable across LP backends, so we
make no recommendation there.

\subsection{Column generation behaviour}
\label{sec:cgbehaviour}

The runtime comparison above establishes which formulation is faster but not why. This section
reports the quantities that explain it: how many iterations each formulation takes, how many columns it
prices in and keeps, and where its wall-clock time goes. \autoref{tbl:cgbehaviour} collects them
per instance family.

Two conventions apply to the table. The counters of a run stopped at the time limit
are truncated (Austin under the path-based formulation reports the 350\,000 columns it managed to
generate in 8 iterations, not the number it would have needed), so the table is restricted to
the 201 (instance, LP backend) cells that \emph{both} formulations solved; this discards 17 of
the 45 transportation path runs, each a run the path-based formulation did not solve, so the restriction understates the effect reported below. And the four column counters
measure different things and do not reconcile into a single balance; the caption states each
one. The same counters are reported per run rather than pooled, including the runs this
section's aggregates exclude, in \supplement.

\begin{table}[htbp]
\centering
\small
\begin{adjustbox}{max width=\textwidth}
\begin{tabular}{llrrrrrrrrrr}
\toprule
 & & & \multicolumn{5}{c}{Columns} & Lazy & \multicolumn{3}{c}{Share of runtime (\%)} \\
\cmidrule(lr){4-8}\cmidrule(lr){10-12}
Family & Form. & Iter. & seeded & generated & purged & final & \emph{of which slack} & rows & master & pricing & sep. \\
\midrule
Grid & Tree & 18 & 466 & 1\,599 & 319 & 1\,500 & 234 & 153 & 97.1 & 1.4 & 0.3 \\
 & Path & 11 & 1\,225 & 2\,123 & 96 & 3\,023 & 214 & 161 & 99.0 & 0.4 & 0.1 \\
\addlinespace
Planar & Tree & 25 & 428 & 3\,789 & 1\,604 & 1\,936 & 214 & 194 & 98.8 & 0.3 & 0.1 \\
 & Path & 11 & 2\,440 & 8\,475 & 525 & 8\,884 & 186 & 203 & 97.4 & 0.2 & 0.8 \\
\addlinespace
Transportation & Tree & 34 & 829 & 8\,196 & 4\,859 & 2\,953 & 415 & 357 & 87.7 & 5.3 & 1.0 \\
 & Path & 14 & 63\,581 & 73\,088 & 7\,446 & 95\,828 & 191 & 387 & 98.9 & 0.4 & 0.2 \\
\addlinespace
Intermodal & Tree & 13 & 31\,904 & 1\,437 & 0 & 33\,611 & 15\,952 & 136 & 22.8 & 65.0 & 9.6 \\
 & Path & 13 & 31\,904 & 1\,446 & 0 & 33\,627 & 15\,952 & 135 & 27.1 & 61.4 & 9.1 \\
\bottomrule
\end{tabular}
\end{adjustbox}
\vspace{1em}
\caption{Column generation behaviour per instance family, over the 201 (instance, LP backend) cells solved by both formulations. Column counts are geometric means over those cells, taken with a shift of one so that a zero count is admissible, since no column is ever purged on the intermodal family, for which an unshifted geometric mean is undefined; the runtime shares are pooled over them. \emph{Seeded} and \emph{final} are the width of the master problem at the first and last iteration, both counting slack columns, and \emph{generated} and \emph{purged} count the structural columns priced in and removed over the whole run. Within a single run the four balance, $\text{final} = \text{seeded} + \text{generated} - \text{purged}$, exactly in 319 of the 420 solved cells and to within 1\% in 400 of them. That identity is a property of each run and not of this summary: the entries below are geometric means taken independently over the cells of a family, and a geometric mean is not additive, so the columns here are not expected to reproduce it. \emph{Of which slack} is the slack part of \emph{final}: one column per structural row when those are fewer than the edges, and otherwise one per lazily added capacity row, as \autoref{sec:colgen} describes. It is a substantial share of the master on the intermodal family, where the seed is one column and one slack per commodity. \emph{Lazy rows} counts the capacity constraints added on demand by the mechanism of \autoref{sec:colgen}. The runtime shares do not sum to 100: the remainder is column management and bookkeeping, which is not attributable to a component.}
\label{tbl:cgbehaviour}
\end{table}

\paragraph{Iterations.} The path-based formulation converges in fewer iterations on every family
where the two differ: 11 against 18 on grid, 11 against 25 on planar, and 14 against 34 on
transportation, as shifted geometric means over the cells both solve; the shift is defined in the
caption of \autoref{tbl:cgbehaviour}. This is the behaviour reported by Jones
et al.~\cite{jones1993multicommodity}, and the reason they gave for preferring the path-based
formulation. Our measurements confirm it. What has changed is that the iteration count no longer
determines the runtime.

\paragraph{Runtime decomposition.} On the grid, planar and transportation
families the restricted master problem accounts for between 87\% and 99\% of wall-clock time under
both formulations, pooling each family's cells rather than averaging over runs; pricing never
exceeds 5.3\% and separation never exceeds 1.1\%. The intermodal family is the sole exception, and
it is reversed: on the same pooled basis 25\% master against 63\% pricing, with a further 9\% in
separation. This measurement is the basis for the per-family strategy choice of
\autoref{sec:colgen}, which selects master-heavy operation on the first three families and
pricing-heavy operation on intermodal; the classification is thus based on measurement rather
than tuning. It also bounds the possible effect of the pricing accelerations on those three
families: with pricing under 5.3\% of runtime, removing its cost entirely would leave the
comparison unchanged.

\paragraph{Master problem size.} Since the master dominates and the path-based formulation needs
fewer master solves, it would be expected to be faster, as it was at the scale of Jones et al.
It is slower here because its master problems are far larger. The path-based master carries one demand constraint per
commodity and is seeded with one column per commodity; the tree-based master carries one demand
constraint per source and is seeded with one column per source. On the transportation family, the
family with the most source sharing, that is a seed of 63\,581 columns against 829 and a final
master of 95\,828 columns against 2\,953, a factor of 32. The difference in master size is far larger
than the difference in pricing activity that produces it: over the same cells the path-based formulation
generates about nine times as many columns as the tree-based one, but ends with a master
thirty-two times larger. Purging does not close the gap, because it cannot remove the
per-commodity seed the path-based master needs to stay feasible. The floor on the path-based master is $|K|$; the floor on the tree-based master is $|S|$.
The tree-based formulation therefore performs a larger number of much cheaper master solves in
place of a smaller number of expensive ones, and at this scale the exchange is favourable.

Austin illustrates the mechanism. Under COPT in CPU mode, with $|K| = 1\,081\,717$ commodities
sharing $|S| = 1\,117$ sources:

\begin{center}
\begin{adjustbox}{max width=\textwidth}
\begin{tabular}{lrrrrrr}
\toprule
 & Iterations & Seeded & Generated & Purged & Final master & \emph{of which slack} \\
\midrule
Path & 8 & 1\,082\,300 & 350\,000 & 17\,349 & 1\,415\,237 & 869 \\
Tree & 223 & 2\,234 & 228\,464 & 206\,319 & 24\,379 & 1\,117 \\
\bottomrule
\end{tabular}
\end{adjustbox}
\end{center}

Austin also shows the two slack regimes side by side: the tree-based master carries one slack
per convexity row, 1\,117 of them, placed before the loop starts, while the path-based master has
more demand rows than edges, so its slacks attach to the capacity constraints instead and arrive
with them, 869 by the time the run stops (the caption of \autoref{tbl:cgbehaviour} states the
rule). On every large transportation instance the rule splits the same way, so on this family slack handling varies
together with the formulation; the mechanism described above does not depend on it, since the path-based master carries one demand row per commodity whatever the slacks
attach to.

The path-based run did not finish: it spent 7\,864 of its 7\,883 seconds inside the LP solver and
completed 8 iterations, because its first master problem already has over a million columns and a
million demand rows before pricing has run once. The tree-based run solved the instance in 559
seconds, of which 492 were master solves spread over 223 iterations (about 2.2 seconds each,
against roughly 980 seconds each for the path-based formulation). It generated 228\,464 columns,
the same order as the 350\,000 the path-based run managed before hitting the limit, but purged
206\,319 of them and finished with a master 58 times smaller. On this instance the tree-based formulation is therefore faster not because it generates fewer
columns (it generates the same order) but because it does not retain them. Over the
transportation family as a whole both effects are present.

The per-iteration column bounds asserted in \autoref{sec:colgen} hold throughout: across all 201
cells no tree-based iteration generated more than $|S|$ columns and no path-based iteration more
than $|K|$.

\paragraph{The intermodal family as a control.} On the intermodal family, the control of
\autoref{sec:instances}, the table shows the two formulations behaving identically, as the same
model must: 13 iterations each, the same seed, and the same final master to within 0.2\%.

The only counter that differs is the time split, which differs by up to 4.3
percentage points, the path-based runs spending 27.1\% of their time in the master against
22.8\%. This is not a formulation effect: the two runs solve the same model but not with
the same pricer, since the stopping bound of \autoref{sec:colgen} is weaker under the path-based
formulation than under the tree-based one (\autoref{sec:ablation}) and saves a median 1.7\% per
pricing problem there against 2.3\%.

\subsection{The bounded pricing criterion}
\label{sec:ablation}

This section measures the one acceleration that has no published counterpart, the bounded
single-source stopping criterion of \autoref{sec:colgen}, in a paired
A/B: the same instances, formulations and backends, run with the criterion switched off and on,
both arms of every comparison in a single session on a single build, because wall clock on
these instances drifts by a few percent between sessions even on cells the criterion provably
leaves unchanged. Both arms of every cell are printed in \supplement. The possible gain from the criterion is
the pricing share times the saving per pricing problem, less the cost in master time of a
shifted trajectory: the criterion returns the same columns but not always along the same iteration
sequence, for the two reasons \autoref{sec:colgen} gives, and where the iteration count changes,
that change dominates the wall-clock difference. The two factors of the gain are anticorrelated
across our families, since the saving grows with the number of commodities per source and the
families with many commodities per source spend almost no time pricing, so we vary one factor
per round: round (a) varies the family, which sets the pricing share, and round (b) the LP
backend, which sets the master share.

\paragraph{Round (a): instance families.}
\autoref{tbl:ablationfamilies} reports 74 cells over the four families, 444 runs at three
repetitions per arm, COPT in GPU mode throughout with the intermodal family additionally run in
CPU mode. The decisive case is grid under the tree-based formulation: the criterion removes more than a
fifth of the pricing time, $-22.4\%$, but wall-clock time rises by $+1.95\%$, because pricing is 3.2\% of it and the
effect is smaller than the session-internal noise of the master term. Planar gives the same null
result, and transportation has the same pricing share of 3.2\% family-wide; the round ran only the 6 of its 9 instances cheap enough to repeat, which
price for 24.3\% of their own runtime against 1.8\% for the three excluded, so its
transportation row is not a family figure, and the grid row, at the same pricing share, stands
in for it. Only intermodal has a large ceiling, at an 80.5\% pricing share in GPU mode and 85.1\% in CPU
mode: wall-clock time falls by 3.6\% and 3.7\%, the gain is entirely in pricing, and it is the
only family whose per-pricing-problem saving, $-2.4\%$ and $-2.6\%$, exceeds its own repetition
noise (spreads of
1.1 and 1.3 percentage points). The table also separates the two bounds. Under the tree-based
formulation the criterion removes 22.4\% of grid's pricing time and 8.9\% of planar's; under
the path-based formulation, on the same instances, it moves pricing by $+0.3\%$ and $-0.8\%$,
because the tree bound tightens each time a sink is settled, whereas the path bound is a
maximum over the unsettled sinks and changes only when the sink carrying the maximum is settled.

\begin{table}[htbp]
\centering
\small
\begin{adjustbox}{max width=\textwidth}
\begin{tabular}{llrrrrrrr}
\toprule
 & & & & \multicolumn{3}{c}{Change, bound on against off (\%)} & \multicolumn{2}{c}{Per pricing problem} \\
\cmidrule(lr){5-7}\cmidrule(lr){8-9}
Family & Form. & Cells & Pricing share & fires & pricing & wall clock & change & noise floor \\
\midrule
Planar & Tree & 9 & 1.6 & 23 & $-8.9$ & $-1.50$ & $-8.0$ (6) & 20.8 \\
 & Path & 9 & 1.0 & 47 & $-0.8$ & $+0.09$ & $0.0$ (8) & 5.0 \\
\addlinespace
Grid & Tree & 15 & 3.2 & 61 & $-22.4$ & $+1.95$ & $0.0$ (6) & 7.1 \\
 & Path & 15 & 1.7 & 64 & $+0.3$ & $+0.10$ & $+2.1$ (9) & 6.3 \\
\addlinespace
Transportation & Tree & 6 & 26.3 & 35 & $-8.7$ & $-1.87$ & $-11.3$ (1) & 2.0 \\
\addlinespace
Intermodal & Tree & 10 & 80.5 & 72 & $-4.5$ & $-3.58$ & $-2.6$ (4) & 1.3 \\
 & Tree, CPU mode & 10 & 85.1 & 72 & $-4.5$ & $-3.74$ & $-2.4$ (4) & 1.1 \\
\bottomrule
\end{tabular}
\end{adjustbox}
\vspace{1em}
\caption{Round (a) of the bounded-pricing ablation: the criterion of \autoref{sec:colgen} switched off and on over four instance families at one LP backend, COPT in GPU mode, three repetitions per arm with both arms run in a single session. A cell is one instance under one formulation, with each arm reduced to the median over its repetitions; the intermodal family is additionally run under COPT in CPU mode, which is the same library at a different barrier execution and so a control on the timing rather than a second result. \emph{Pricing share} is the baseline arm's pricing time as a percentage of its wall clock, and it is the ceiling on anything the criterion can convert into runtime; \emph{fires} is the median percentage of pricing problems on which the stopping criterion triggered. The two \emph{change} columns are ratios of summed seconds over a group's cells, not means of per-cell percentages, so each is weighted as a wall-clock total would be. \emph{Per pricing problem} is the same measurement divided by the number of pricing problems each arm actually solved, which unlike a total does not move when the two arms price a different number of sources; it is reported over the cells on which the two arms ran the same trajectory (equal iteration and column counts), with that count in parentheses, and \emph{noise floor} is the baseline arm's own repetition-to-repetition spread on it. Where that floor is of the same order as the change, as it is on grid and planar, the estimate is reported but carries no information. Both arms agree on the LP optimum on every cell, to a worst relative difference of $6.95 \times 10^{-5}$, inside the $\epsilon_{\mathrm{rel}} = 10^{-4}$ the runs are solved to.}
\label{tbl:ablationfamilies}
\end{table}

\paragraph{Round (b): LP backends.}
Round (a) used the fastest LP available, which minimizes the cost of a shifted trajectory, and
that term is the one that could make the result an artefact of the solver. Round (b),
\autoref{tbl:ablationbackends}, takes the intermodal family and both formulations across all
five LP configurations: 100 cells, 240 runs, no failures. Pricing time falls under every
backend, in all ten solver--formulation groups, by between 1.3\% and 6.1\%; that is a property
of the pricer, independent of which code solves the master. Wall-clock time follows
only where the master is not the bottleneck: $-5.4\%$ and $-6.0\%$ under COPT and $-2.2\%$
under MOSEK, on master shares of 3.5\%, 9.5\% and 6.9\%, against $+2.2\%$ under cuOpt and
$+15.7\%$ under HiGHS, on 13.8\% and 51.8\%. The two positive entries are due to the trajectory term,
not to a cost of the criterion: on the 42 cells whose trajectory was unchanged, master time is flat to
within a percentage point under every backend, pricing falls 1.4\% to 2.5\%, and wall clock
falls 0.8\% to 2.2\%, negative under all five; over all 100 cells the criterion is faster on 76,
at a median of $-1.7\%$. The trajectory term is largest under HiGHS, whose per-iteration master cost
is the highest: BUS-23688-0 goes from 20 iterations to 44 under the path-based formulation
and from 38 to 20 under the tree-based one, for $+135.5\%$ and $-48.7\%$ on an identical
column set, and repeating the round in a second session reversed the sign of both HiGHS
totals.\footnote{This is the one statement in this section not derivable from the
data accompanying this paper, which holds a single measurement per round. The second session's
run logs are not part of the archived artefact, which ships the derived summaries only; we report
the reversal because a single measurement cannot establish its own stability, and it bears
directly on how much weight the HiGHS column of \autoref{tbl:ablationbackends} can carry.} A wall-clock total over such cells is not a stable quantity, and we
do not report one as the criterion's effect on that backend. The per-pricing-problem saving
agrees with round (a), a median of $-2.0\%$ over the 42 stable cells against round (a)'s
$-2.4\%$. On a family that spends 76\% to 84\% of its runtime in pricing under four of the five
backends, the gain justifies enabling the criterion; on the other three families it does not, and
the criterion is enabled on intermodal alone.

\begin{table}[htbp]
\centering
\small
\begin{adjustbox}{max width=\textwidth}
\begin{tabular}{lrrrrrrrrr}
\toprule
 & \multicolumn{6}{c}{All cells} & \multicolumn{3}{c}{Trajectory unchanged} \\
\cmidrule(lr){2-7}\cmidrule(lr){8-10}
 & \multicolumn{2}{c}{Share of runtime (\%)} & \multicolumn{2}{c}{Change (\%)} & & & & \multicolumn{2}{c}{Change (\%)} \\
\cmidrule(lr){2-3}\cmidrule(lr){4-5}\cmidrule(lr){9-10}
LP backend & pricing & master & pricing & wall clock & faster & cells & cells & pricing & wall clock \\
\midrule
HiGHS & 42.1 & 51.8 & $-2.5$ & $+15.7$ & 14 & 20 & 8 & $-1.4$ & $-0.8$ \\
MOSEK & 82.9 & 6.9 & $-3.8$ & $-2.2$ & 16 & 20 & 8 & $-2.5$ & $-2.2$ \\
cuOpt & 76.1 & 13.8 & $-2.7$ & $+2.2$ & 14 & 20 & 10 & $-1.7$ & $-1.5$ \\
COPT (CPU) & 84.2 & 3.5 & $-5.1$ & $-5.4$ & 16 & 20 & 8 & $-2.1$ & $-1.9$ \\
COPT (GPU) & 78.6 & 9.5 & $-4.9$ & $-6.0$ & 16 & 20 & 8 & $-1.8$ & $-1.6$ \\
\midrule
All & 68.3 & 22.2 & $-3.8$ & $+3.1$ & 76 & 100 & 42 & $-1.9$ & $-1.6$ \\
\bottomrule
\end{tabular}
\end{adjustbox}
\vspace{1em}
\caption{Round (b) of the bounded-pricing ablation: the one family round (a) left standing, run under all five LP backends, both formulations, with both arms of every comparison in a single session on a single build. Columns are as in \autoref{tbl:ablationfamilies}. \emph{Faster} counts the cells on which the bound reduced wall clock, out of the \emph{cells} beside it. The right-hand panel restricts to the cells on which the two arms ran the same trajectory, which removes the $\pm$iteration term and leaves the pricing effect alone: there the master time is flat to within a percentage point under every backend and the pricing saving passes through to the clock. The pricing saving in the left-hand panel is negative without exception, which is the result that does not depend on the backend. The wall-clock column beside it does, and its two positive entries are not a cost the pricer imposes: on cuOpt and HiGHS the iteration count moves in both directions and the master problem is expensive enough that it decides the total, as \autoref{sec:ablation} explains. The \emph{All} row is included for completeness and should not be read as a summary: it is a sum of seconds across backends that differ by a factor of two in speed, so the slowest dominates it.}
\label{tbl:ablationbackends}
\end{table}

\subsection{The planar2500 counterexample}
\label{sec:planar2500}

The mechanism of \autoref{sec:cgbehaviour}, more iterations but on a much smaller and therefore
much cheaper master, accounts for the aggregate result and for the transportation family in
particular. planar2500 is the instance on which it does not hold: the tree-based master there is seventeen
times smaller than the path-based one and is nonetheless up to eight times more expensive to
solve. As the only clear counterexample in the study, it is analysed in detail here.

\begin{figure}[htbp]
  \centering
  \includegraphics[width=\textwidth]{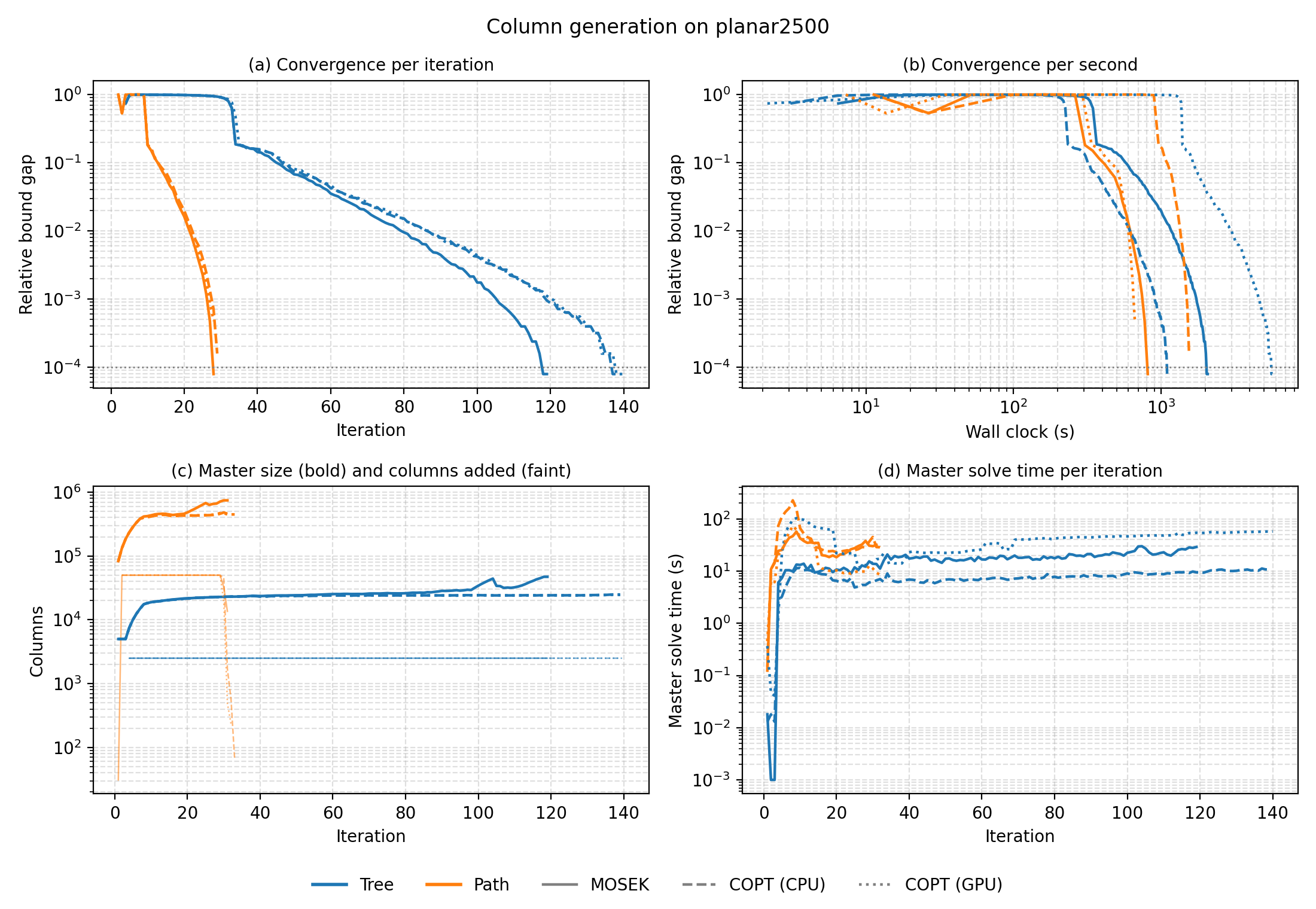}
  \caption{Column generation on planar2500 under the three backends that solve it under both formulations. (a) and (b) plot the relative bound gap $(z_{\mathrm{RMP}} - z_{\mathrm{LB}})/|z_{\mathrm{RMP}}|$ against iteration and against wall clock; the opening iterations are masked, where lazy separation leaves the master a relaxation and the difference is not a gap. (c) gives the master size (bold) and the columns added per iteration (faint). (d) gives the master solve time per iteration}
  \label{fig:convergence-planar2500}
\end{figure}

\autoref{fig:convergence-planar2500} traces both formulations iteration by iteration under the
three backends that solve the instance under both, and \autoref{tbl:planar2500} collects the
exact counters. Under cuOpt the path-based formulation solves it in 544 seconds and the tree-based
one does not finish; under HiGHS neither does.

\begin{table}[htbp]
\centering
\caption{planar2500 under the three backends that solve it under both formulations. Master columns are averaged over the run, and the cost per column is the master solve time per iteration divided by that average}
\label{tbl:planar2500}
\small
\begin{adjustbox}{max width=\textwidth}
\begin{tabular}{llrrrrr}
\toprule
Backend & Form. & Iterations & Master cols. & Cost/col. ($\mu$s) & Master (s) & Total (s) \\
\midrule
MOSEK & Tree & 119 & 25\,797 & 666 & 2\,044 & 2\,071 \\
 & Path & 31 & 469\,684 & 61 & 885 & 909 \\
COPT (CPU) & Tree & 139 & 22\,759 & 343 & 1\,085 & 1\,122 \\
 & Path & 33 & 392\,729 & 127 & 1\,643 & 1\,683 \\
COPT (GPU) & Tree & 140 & 22\,730 & 1\,786 & 5\,684 & 5\,722 \\
 & Path & 32 & 391\,072 & 54 & 672 & 712 \\
\bottomrule
\end{tabular}
\end{adjustbox}
\end{table}

\paragraph{Decomposition of the runtime ratio.} The master problem accounts for between 94\% and 99\%
of wall-clock time in all six runs, so the runtime ratio is the master time ratio, and that ratio
factors into three terms. Writing $n$ for the iteration count, $m$ for the mean master size and
$c$ for what one column costs the barrier in one solve, the master time of a run is $n\,m\,c$ and

\[
\frac{T_{\text{tree}}}{T_{\text{path}}}
 = \underbrace{\frac{n_{\text{tree}}}{n_{\text{path}}}}_{\text{iterations}}
 \times \underbrace{\frac{m_{\text{tree}}}{m_{\text{path}}}}_{\text{master size}}
 \times \underbrace{\frac{c_{\text{tree}}}{c_{\text{path}}}}_{\text{cost per column}}.
\]

Two of the three factors are stable across backends, and these are the two determined by the
formulation. The tree-based formulation takes 3.8 to 4.4 times as many iterations, and its master
holds 5.5\% to 5.8\% of the columns. Multiplied out, those two would put the tree-based
formulation ahead by a factor of about four. The third factor determines the outcome, and it
is not stable: a tree column costs the barrier 2.7 times what a path column costs under COPT
in CPU mode, 11 times under MOSEK and 33 times under COPT in GPU mode. The products are 0.66, 2.31
and 8.45, which are the observed master time ratios. On planar2500, therefore, the
tree-based formulation does not converge badly (panel (a) shows the same convergence behaviour as
on the other instances); its columns are expensive.

\paragraph{Column density.} The cost of a column in a barrier solve is determined by its number
of non-zeros. A path column carries one demand entry and one capacity entry
per edge of a single path; a tree column carries them per edge of a shortest-path tree, and
planar2500 is the extreme case of that difference: every vertex is a source,
$|S| = |V| = 2\,500$, so a tree rooted at a source spans up to $|V| - 1 = 2\,499$ of the
network's 12\,990 edges. Lazy separation keeps the number of capacity rows small in both masters
(under COPT in CPU mode they end with 1\,777 and 1\,791 of them), so a single tree column can
have a non-zero in more than 40\% of the 4\,277 rows of the tree-based master, while a path
column touches a handful of the 83\,221 rows of the path-based one. The tree-based master is the
smaller matrix and the denser one, which is the mechanism Jones et al.\ identified when they
concluded that ``the larger and sparser (ODP) formulation is easier to solve than the smaller
and denser (DSP) formulation in a decomposition framework''
\cite[p.~112]{jones1993multicommodity}. planar2500 allows it to be quantified.

\begin{figure}[htbp]
  \centering
  \includegraphics[width=0.85\textwidth]{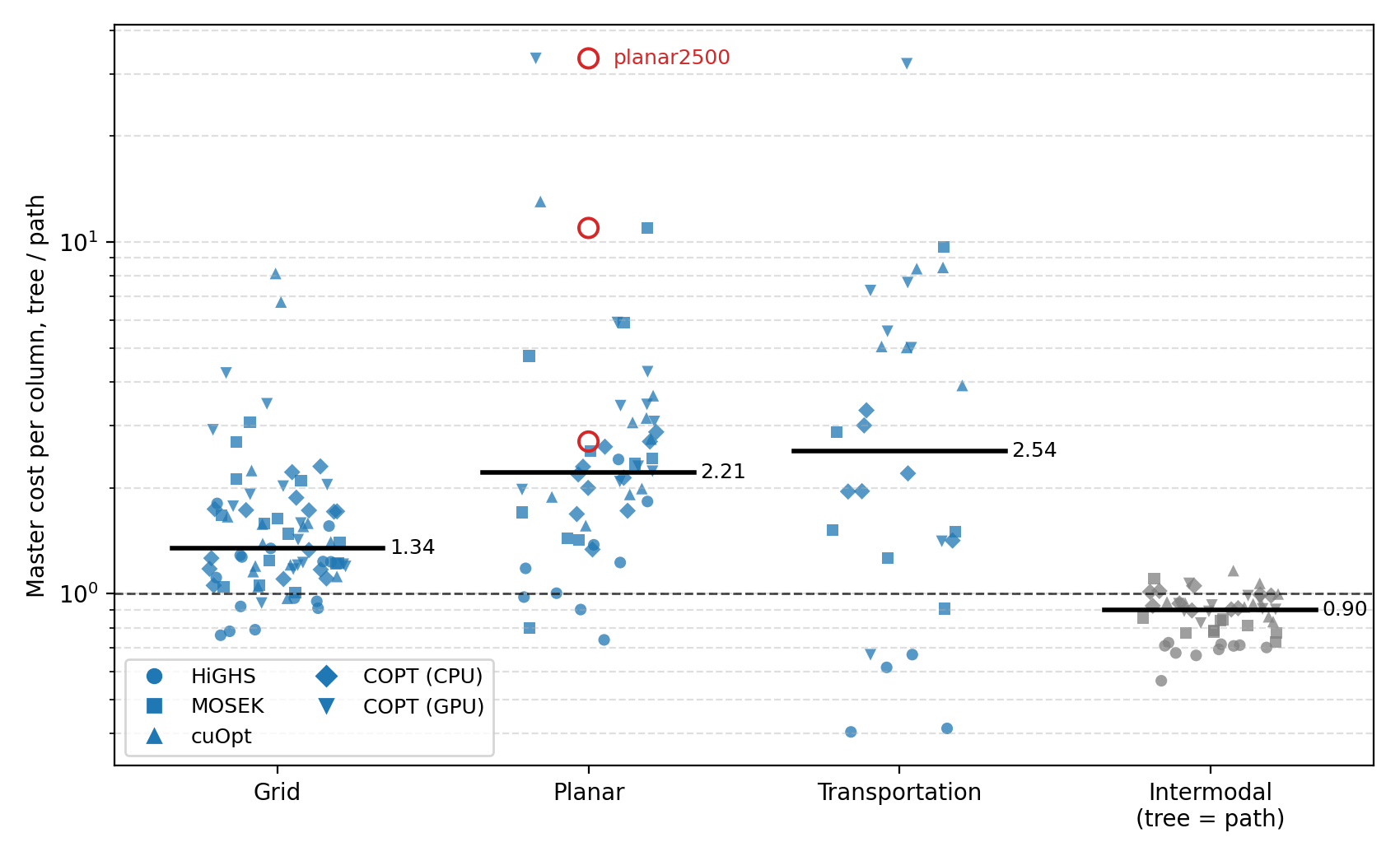}
  \caption{Master solve cost per column, tree-based over path-based, for each (instance, LP backend) cell that both formulations solved; bars are family medians. On the intermodal family a shortest-path tree is a path, so the two formulations build the same columns}
  \label{fig:master-cost-premium}
\end{figure}

\autoref{fig:master-cost-premium} tests this interpretation against the whole study, over the 201
(instance, LP backend) cells both formulations solved. The premium is above 1 in 132 of the 151
cells where the tree-based columns are genuinely trees, with family medians of 1.34 on grid, 2.21
on planar and 2.54 on transportation. On the intermodal family, the control of
\autoref{sec:instances}, it is above 1 in 7 of 50 cells and its median is 0.90. This contrast is the closest approximation to a controlled test the study affords: the premium
appears where the two formulations build different columns and disappears where they build the same
ones. Two biases in the measurement both work against this interpretation. Barrier cost
is superlinear in problem size, so dividing by the column count favours the larger master, which
is the path-based one in every cell; and restricting to cells both formulations solved discards 17
transportation path-based runs that did not finish, every one of them a case where the path-based
master had become unmanageable. The premium also rises with the number of
commodities per source, which is what sets how far a single tree must spread (Spearman 0.68 over
201 cells), though that quantity is the reciprocal of the source ratio of
\autoref{fig:heatmap-path-tree} and hence a size proxy, not an independent variable. The cells are
five backends applied to the same instances rather than independent samples, so the coefficient
describes this suite and is not offered as an inferential statistic.

\paragraph{Open questions.} Column density explains why the premium exists; it does not
explain how large it is on any given run. The clearest evidence comes from planar2500 itself: COPT
in CPU mode and in GPU mode solve the same sequence of master problems, and the premium is 2.7 in
one and 33 in the other, a factor of twelve from the barrier's execution mode alone. Which property of the barrier is responsible is a question about the LP solvers rather than the
formulations, and our instrumentation does not answer it: we record master solve times but not per-column non-zero
counts, so density is the best-supported candidate rather than a confirmed mechanism. Two points are established. The premium is not a single-solver artefact, being 2.7, 11 and 33 under the three
backends that solve the instance, although only under two of them is it large enough to decide the run. And the counterexample is confined to runtime: the tree-based formulation still uses less
memory here, 1.9 to 3.0\,GB against 4.2 to 5.0\,GB, under all three.

\section{Conclusion}
\label{sec:conclusion}

We have revisited the tree-based formulation for the multi-commodity flow problem of Jones et
al.~\cite{jones1993multicommodity}, which aggregates flows by shared source and represents them
as trees, so that the master problem has $O(|S|)$ rather than $O(|K|)$ demand constraints, where
$|S| \ll |K|$ in many practical applications. Re-evaluating it at modern scale, we find:

\begin{itemize}
  \item On large instances with many shared sources ($|S| \ll |K|$) the tree-based formulation
        scales substantially better than path-based column generation: between 1.42 and 1.93
        times faster on the shifted geometric mean over 44 instances depending on the LP
        backend, reaching a factor of 99 on the largest transportation instance both
        formulations solve, with the path-based formulation failing outright on
        17 of the 45 transportation runs.
  \item The ordering does not depend on the interior point code: four barrier implementations
        in five configurations, open source and commercial, CPU and GPU, reproduce it.
  \item Both decompositions substantially outperform direct LP solving, which is between 15 and
        28 times slower on the mean and, on nine of the 44 instances, never gets as far as a
        model to solve.
  \item The mechanism is measurable. The path-based formulation still converges in fewer
        iterations, as Jones et al.\ observed, but the master problem accounts for 87\% to 99\%
        of runtime on three of the four families, and the path-based master ends 32 times
        larger than the tree-based one on the transportation family. A larger number of cheap master solves is faster than a smaller number of expensive
        ones.
  \item Accelerating the pricer is effective only where pricing is the bottleneck. Our stopping
        criterion reduces pricing time under every LP backend without exception, and on three
        of the four families this has no effect on runtime, because pricing is at most 5.3\% of
        their runtime. On the family where pricing is 76\% to 84\% of the runtime the criterion
        saves 2\% to 6\% under three of the four backends in that band, while under cuOpt a
        shifted iteration trajectory leaves the total 2\% higher; restricted to the cells where the
        trajectory held, it is negative under all five.
  \item The memory advantage is smaller than the runtime advantage, a geometric mean of
        1.09 where both formulations solve, rising to 1.37 on the transportation family, and it
        becomes decisive only where the path-based master no longer fits: Sydney under MOSEK
        was terminated by the kernel at 95.8\,GB, while the tree-based formulation solved the
        same instance in 76 seconds using 11.2\,GB.
  \item The tree-based formulation is not uniformly better. On planar2500, where $|S| = |V|$,
        a shortest-path tree spans a fifth of the arcs and fills over 40\% of the tree master's
        rows; its columns cost up to 33 times more per column to process, and the path-based
        formulation wins by up to a factor of eight. This is the mechanism Jones et al.\
        identified, at a scale where it can be quantified, and it is the second condition our
        rule needs: source sharing is necessary but not sufficient, and the advantage disappears when
        $|S|$ approaches $|V|$.
\end{itemize}

These results are subject to three limitations. The suite is 44 instances in four families, and the
reversal is carried by the nine transportation instances, the only ones combining a low source
ratio with $|K|$ reaching into the millions; the grid and planar families sit near
$|S| \approx |K|$ and the two formulations are within a small factor of each other there, and the intermodal family is
a control on which the two formulations are the same model. Every backend is an interior point
code with crossover disabled, a simplex master being impractical at this size
(\autoref{sec:performance}). And the two formulations are
compared under the slack regime each selects, which on the transportation family means the
comparison varies slack handling alongside the formulation, though not the row floor that drives
the result.

These results reverse the practical conclusion of Jones et al.\ at modern scale. Their
convergence finding is confirmed: the path-based formulation still
converges in fewer master iterations, here by a factor of two to three. What is reversed is the
wall-clock consequence they drew from it, and the reversal is not explained by the degree of
source aggregation. Their instances carry between 3 and 79 commodities per source, a range that
spans and exceeds our grid and planar families, yet the path-based formulation won there by its
widest margins at the highest sharing ratios, the opposite of what a sharing-driven account
would predict. The two studies differ in absolute scale. Their largest restricted
master problem held 16{,}431 columns and 5{,}808 rows, a size at which the master solve is a
small share of the total and the runtime tracks the iteration count, which the path-based
formulation minimizes. Our path-based masters reach 1.4 million columns, and once the master
solve dominates the clock, the formulation that keeps the master small wins despite converging
more slowly.

Scale is one of two factors that have changed since 1993, and the one whose effect we measure
directly. The other is the LP technology applied to the master: Jones et al.\ solved theirs with
CPLEX, at the time a simplex code, whose warm start across column generation iterations lowers
the per-iteration master cost and
so shifts weight from master size towards the iteration count, the axis on which the path-based
formulation wins. At their scale the simplex was the appropriate choice; at ours it is not an option, for the
reasons given in \autoref{sec:performance}. The finding is therefore an empirical one about the
problem sizes reachable today and the master technology that reaches them: the ordering Jones
et al.\ reported is reversed, and we do not claim their conclusion was wrong for the instances
and the solvers they had.

Several directions for further work remain. The pricing problem could benefit from
edge filtering, path generation heuristics, and warm-starting techniques, but the master problem
is the bottleneck on three of the four instance families and the primary target: first-order
methods and further GPU acceleration are the candidates, and the spread we observe
between LP backends suggests there is still room there. Extending the tree-based decomposition
to integer multi-commodity flow, where the master problem sits inside a branch-and-price search
and its size compounds across nodes, is the direction we consider most promising.

\bibliographystyle{unsrturl}
\bibliography{references}

\begin{thebibliography}{10}

\bibitem{jones1993multicommodity}
Kim~L. Jones, Irvin~J. Lustig, Judith~M. Farvolden, and Warren~B. Powell.
\newblock Multicommodity network flows: {T}he impact of formulation on
  decomposition.
\newblock {\em Mathematical Programming}, 62:95--117, 1993.
\newblock \href {https://doi.org/10.1007/BF01585162}
  {\path{doi:10.1007/BF01585162}}.

\bibitem{ahuja1993network}
Ravindra~K. Ahuja, Thomas~L. Magnanti, and James~B. Orlin.
\newblock {\em Network Flows: Theory, Algorithms, and Applications}.
\newblock Prentice Hall, Englewood Cliffs, NJ, 1993.

\bibitem{salimifard2022mcnf}
Khodakaram Salimifard and Sara Bigharaz.
\newblock The multicommodity network flow problem: state of the art
  classification, applications, and solution methods.
\newblock {\em Operational Research}, 22(1):1--47, 2022.
\newblock \href {https://doi.org/10.1007/s12351-020-00564-8}
  {\path{doi:10.1007/s12351-020-00564-8}}.

\bibitem{dantzig1960decomposition}
George~B. Dantzig and Philip Wolfe.
\newblock Decomposition principle for linear programs.
\newblock {\em Operations Research}, 8(1):101--111, 1960.
\newblock \href {https://doi.org/10.1287/opre.8.1.101}
  {\path{doi:10.1287/opre.8.1.101}}.

\bibitem{desrosiers2024branch}
Jacques Desrosiers, Marco L{\"u}bbecke, Guy Desaulniers, and Jean~Bertrand
  Gauthier.
\newblock {\em Branch-and-Price}.
\newblock Springer, Cham, 2026.
\newblock \href {https://doi.org/10.1007/978-3-031-96917-1}
  {\path{doi:10.1007/978-3-031-96917-1}}.

\bibitem{uchoa2024column}
Eduardo Uchoa, Artur Pessoa, and Lorenza Moreno.
\newblock {Optimizing with Column Generation}: Advanced branch-cut-and-price
  algorithms ({Part I}).
\newblock Technical Report L-2024-3, Cadernos do LOGIS-UFF, Universidade
  Federal Fluminense, Engenharia de Produ{\c{c}}{\~a}o, August 2024.

\bibitem{gondzio2016primaldual}
Jacek Gondzio, Pablo Gonz{\'a}lez-Brevis, and Pedro Munari.
\newblock Large-scale optimization with the primal--dual column generation
  method.
\newblock {\em Mathematical Programming Computation}, 8:47--82, 2016.
\newblock \href {https://doi.org/10.1007/s12532-015-0090-6}
  {\path{doi:10.1007/s12532-015-0090-6}}.

\bibitem{babonneau2006active}
Fr{\'e}d{\'e}ric Babonneau, Olivier du~Merle, and Jean-Philippe Vial.
\newblock Solving large-scale linear multicommodity flow problems with an
  active set strategy and {Proximal-ACCPM}.
\newblock {\em Operations Research}, 54(1):184--197, 2006.
\newblock \href {https://doi.org/10.1287/opre.1050.0262}
  {\path{doi:10.1287/opre.1050.0262}}.

\bibitem{zhang2025gpu}
Fangzhao Zhang and Stephen Boyd.
\newblock Solving large multicommodity network flow problems on {GPU}s.
\newblock {\em Automation and Remote Control}, 86(8):782--795, 2025.
\newblock \href {https://arxiv.org/abs/2501.17996} {\path{arXiv:2501.17996}},
  \href {https://doi.org/10.48550/arXiv.2501.17996}
  {\path{doi:10.48550/arXiv.2501.17996}}.

\bibitem{yin2020unified}
Ping Yin, Steven Diamond, Bill Lin, and Stephen Boyd.
\newblock Network optimization for unified packet and circuit switched
  networks.
\newblock {\em Optimization and Engineering}, 21(1):159--180, 2020.
\newblock \href {https://doi.org/10.1007/s11081-019-09439-0}
  {\path{doi:10.1007/s11081-019-09439-0}}.

\bibitem{bargera2002origin}
Hillel Bar-Gera.
\newblock Origin-based algorithm for the traffic assignment problem.
\newblock {\em Transportation Science}, 36(4):398--417, 2002.
\newblock \href {https://doi.org/10.1287/trsc.36.4.398.549}
  {\path{doi:10.1287/trsc.36.4.398.549}}.

\bibitem{lienkamp2024column}
Benedikt Lienkamp and Maximilian Schiffer.
\newblock Column generation for solving large scale multi-commodity flow
  problems for passenger transportation.
\newblock {\em European Journal of Operational Research}, 314(2):703--717,
  2024.
\newblock \href {https://doi.org/10.1016/j.ejor.2023.09.019}
  {\path{doi:10.1016/j.ejor.2023.09.019}}.

\bibitem{mcfcg}
Simon Spoorendonk and Bj{\o}rn Petersen.
\newblock {mcfcg}: Column generation for minimum-cost multicommodity flow,
  2026.
\newblock Software, release v0.1.0 (MIT). The version DOI archives the exact
  release behind the results reported here; concept DOI 10.5281/zenodo.22058613
  resolves to the latest release. Source at
  \url{https://github.com/spoorendonk/mcfcg}.
\newblock \href {https://doi.org/10.5281/zenodo.22058614}
  {\path{doi:10.5281/zenodo.22058614}}.

\bibitem{barnhart1994column}
Cynthia Barnhart, Christopher~A. Hane, Ellis~L. Johnson, and Gabriele
  Sigismondi.
\newblock A column generation and partitioning approach for multi-commodity
  flow problems.
\newblock {\em Telecommunication Systems}, 3:239--258, 1994.
\newblock \href {https://doi.org/10.1007/BF02110307}
  {\path{doi:10.1007/BF02110307}}.

\bibitem{barnhart2000branch}
Cynthia Barnhart, Christopher~A. Hane, and Pamela~H. Vance.
\newblock Using branch-and-price-and-cut to solve origin-destination integer
  multicommodity flow problems.
\newblock {\em Operations Research}, 48(2):318--326, 2000.
\newblock \href {https://doi.org/10.1287/opre.48.2.318.12378}
  {\path{doi:10.1287/opre.48.2.318.12378}}.

\bibitem{ford1958suggested}
Lester~R. Ford, Jr. and Delbert~R. Fulkerson.
\newblock A suggested computation for maximal multi-commodity network flows.
\newblock {\em Management Science}, 5(1):97--101, 1958.
\newblock \href {https://doi.org/10.1287/mnsc.5.1.97}
  {\path{doi:10.1287/mnsc.5.1.97}}.

\bibitem{dijkstra1959note}
Edsger~W. Dijkstra.
\newblock A note on two problems in connexion with graphs.
\newblock {\em Numerische Mathematik}, 1:269--271, 1959.
\newblock \href {https://doi.org/10.1007/BF01386390}
  {\path{doi:10.1007/BF01386390}}.

\bibitem{hart1968formal}
Peter~E. Hart, Nils~J. Nilsson, and Bertram Raphael.
\newblock A formal basis for the heuristic determination of minimum cost paths.
\newblock {\em {IEEE} Transactions on Systems Science and Cybernetics},
  4(2):100--107, 1968.
\newblock \href {https://doi.org/10.1109/TSSC.1968.300136}
  {\path{doi:10.1109/TSSC.1968.300136}}.

\bibitem{goldberg2005computing}
Andrew~V. Goldberg and Chris Harrelson.
\newblock Computing the shortest path: {A*} search meets graph theory.
\newblock In {\em Proceedings of the 16th Annual {ACM-SIAM} Symposium on
  Discrete Algorithms ({SODA})}, pages 156--165, 2005.

\bibitem{huangfu2018parallelizing}
Qi~Huangfu and J.~A.~J. Hall.
\newblock Parallelizing the dual revised simplex method.
\newblock {\em Mathematical Programming Computation}, 10(1):119--142, 2018.
\newblock \href {https://doi.org/10.1007/s12532-017-0130-5}
  {\path{doi:10.1007/s12532-017-0130-5}}.

\bibitem{zanetti2025factorisation}
Filippo Zanetti and Jacek Gondzio.
\newblock A factorisation-based regularised interior point method using the
  augmented system.
\newblock {\em arXiv preprint}, 2025.
\newblock \href {https://arxiv.org/abs/2508.04370} {\path{arXiv:2508.04370}},
  \href {https://doi.org/10.48550/arXiv.2508.04370}
  {\path{doi:10.48550/arXiv.2508.04370}}.

\bibitem{mosek}
{MOSEK ApS}.
\newblock {\em {MOSEK} Optimizer {API} for {C}. Version 11.0}.
\newblock MOSEK ApS, 2025.
\newblock Accessed 2026-08-11,
  \url{https://docs.mosek.com/11.0/capi/index.html}.

\bibitem{cuopt}
{NVIDIA}.
\newblock {NVIDIA cuOpt}: {GPU} accelerated decision optimization, 2026.
\newblock Open source (Apache-2.0), also hosted as a COIN-OR project. Accessed
  2026-08-11, \url{https://github.com/NVIDIA/cuopt}.

\bibitem{ge2022copt}
Dongdong Ge, Qi~Huangfu, Zizhuo Wang, Jian Wu, and Yinyu Ye.
\newblock Cardinal {O}ptimizer ({COPT}) user guide.
\newblock {\em arXiv preprint}, 2022.
\newblock \href {https://arxiv.org/abs/2208.14314} {\path{arXiv:2208.14314}},
  \href {https://doi.org/10.48550/arXiv.2208.14314}
  {\path{doi:10.48550/arXiv.2208.14314}}.

\bibitem{mcf-instances}
{CommaLAB, University of Pisa}.
\newblock {Multicommodity Flow Problems}, 2025.
\newblock Accessed 2025-08-27,
  \url{https://commalab.di.unipi.it/datasets/mmcf/}.

\bibitem{transportationnetworks2025}
{Transportation Networks for Research Core Team}.
\newblock {Transportation Networks for Research}, 2025.
\newblock GitHub repository. commit \#d1639b4, Accessed 2025-08-27,
  \url{https://github.com/bstabler/TransportationNetworks}.

\bibitem{achterberg2007constraint}
Tobias Achterberg.
\newblock {\em Constraint Integer Programming}.
\newblock PhD thesis, Technische Universit{\"a}t Berlin, 2007.
\newblock \href {https://doi.org/10.14279/depositonce-1634}
  {\path{doi:10.14279/depositonce-1634}}.

\end{thebibliography}

\appendix
\section{Appendix}
\label{sec:appendix}
\subsection{Instance details}
\autoref{tbl:coefs} shows the coefficients used for the generation of min-cost transportation network instances. \autoref{tbl:instancedetails} shows the instance details.

\begin{table}[htbp]
  \centering
\begin{tabular}{lr}
\toprule
Instance & Coefficient \\
\midrule
\multicolumn{2}{l}{\textbf{Transportation}} \\
Austin & 6.0 \\
Barcelona & 5050.0 \\
BerlinCenter & 0.5 \\
Birmingham & 0.9 \\
ChicagoRegional & 4.1 \\
ChicagoSketch & 2.4 \\
Philadelphia & 7.0 \\
Sydney & 1.9 \\
Winnipeg & 2000.0 \\
\addlinespace
\bottomrule
\end{tabular}
\vspace{1em}
  \caption{Coefficients used for the generation of min-cost transportation network instances}
  \label{tbl:coefs}
\end{table}

\begin{table}[htbp]
\centering
\begin{adjustbox}{max width=\textwidth}
\begin{tabular}{lrrrrrrr}
\toprule
Instance & Nodes & Edges & Commodities & Sources & Variables & Constraints & Non-zeros \\
\midrule
\multicolumn{8}{l}{\textbf{Grid}} \\
grid1 & 25 & 80 & 50 & 23 & 1.8e+03 & 6.6e+02 & 5.5e+03\\
grid2 & 25 & 80 & 100 & 25 & 2.0e+03 & 7.0e+02 & 6.0e+03\\
grid3 & 100 & 360 & 50 & 40 & 1.4e+04 & 4.4e+03 & 4.3e+04\\
grid4 & 100 & 360 & 100 & 63 & 2.3e+04 & 6.7e+03 & 6.8e+04\\
grid5 & 225 & 840 & 100 & 83 & 7.0e+04 & 2.0e+04 & 2.1e+05\\
grid6 & 225 & 840 & 200 & 135 & 1.1e+05 & 3.1e+04 & 3.4e+05\\
grid7 & 400 & 1,520 & 400 & 247 & 3.8e+05 & 1.0e+05 & 1.1e+06\\
grid8 & 625 & 2,400 & 500 & 343 & 8.2e+05 & 2.2e+05 & 2.5e+06\\
grid9 & 625 & 2,400 & 1,000 & 495 & 1.2e+06 & 3.1e+05 & 3.6e+06\\
grid10 & 625 & 2,400 & 2,000 & 593 & 1.4e+06 & 3.7e+05 & 4.3e+06\\
grid11 & 625 & 2,400 & 4,000 & 625 & 1.5e+06 & 3.9e+05 & 4.5e+06\\
grid12 & 900 & 3,480 & 6,000 & 899 & 3.1e+06 & 8.1e+05 & 9.4e+06\\
grid13 & 900 & 3,480 & 12,000 & 900 & 3.1e+06 & 8.1e+05 & 9.4e+06\\
grid14 & 1,225 & 4,760 & 16,000 & 1,225 & 5.8e+06 & 1.5e+06 & 1.7e+07\\
grid15 & 1,225 & 4,760 & 32,000 & 1,225 & 5.8e+06 & 1.5e+06 & 1.7e+07\\
\addlinespace
\multicolumn{8}{l}{\textbf{Planar}} \\
planar30 & 30 & 150 & 92 & 29 & 4.4e+03 & 1.0e+03 & 1.3e+04\\
planar50 & 50 & 250 & 267 & 50 & 1.2e+04 & 2.8e+03 & 3.8e+04\\
planar80 & 80 & 440 & 543 & 80 & 3.5e+04 & 6.8e+03 & 1.1e+05\\
planar100 & 100 & 532 & 1,085 & 100 & 5.3e+04 & 1.1e+04 & 1.6e+05\\
planar150 & 150 & 850 & 2,239 & 150 & 1.3e+05 & 2.3e+04 & 3.8e+05\\
planar300 & 300 & 1,680 & 3,584 & 300 & 5.0e+05 & 9.2e+04 & 1.5e+06\\
planar500 & 500 & 2,842 & 3,525 & 500 & 1.4e+06 & 2.5e+05 & 4.3e+06\\
planar800 & 800 & 4,388 & 12,756 & 800 & 3.5e+06 & 6.4e+05 & 1.1e+07\\
planar1000 & 1,000 & 5,200 & 20,026 & 1,000 & 5.2e+06 & 1.0e+06 & 1.6e+07\\
planar2500 & 2,500 & 12,990 & 81,430 & 2,500 & 3.2e+07 & 6.3e+06 & 9.7e+07\\
\addlinespace
\multicolumn{8}{l}{\textbf{Transportation}} \\
Austin & 7,388 & 18,961 & 1,081,717 & 1,117 & 2.1e+07 & 8.3e+06 & 6.4e+07\\
Barcelona & 1,020 & 2,522 & 7,922 & 97 & 2.4e+05 & 1.0e+05 & 7.3e+05\\
BerlinCenter & 12,981 & 28,376 & 49,688 & 865 & 2.5e+07 & 1.1e+07 & 7.4e+07\\
Birmingham & 14,639 & 33,937 & 471,637 & 898 & 3.0e+07 & 1.3e+07 & 9.1e+07\\
ChicagoRegional & 12,982 & 39,018 & 2,297,945 & 1,771 & 6.9e+07 & 2.3e+07 & 2.1e+08\\
ChicagoSketch & 933 & 2,950 & 93,513 & 386 & 1.1e+06 & 3.6e+05 & 3.4e+06\\
Philadelphia & 13,389 & 40,003 & 1,151,166 & 1,489 & 6.0e+07 & 2.0e+07 & 1.8e+08\\
Sydney & 33,113 & 75,379 & 3,343,560 & 3,257 & 2.5e+08 & 1.1e+08 & 7.4e+08\\
Winnipeg & 1,052 & 2,836 & 4,345 & 135 & 3.8e+05 & 1.4e+05 & 1.1e+06\\
\addlinespace
\multicolumn{8}{l}{\textbf{Intermodal}} \\
BUS-2632-0 & 119,865 & 397,362 & 2,628 & 2,628 & 1.0e+09 & 3.2e+08 & 2.6e+09\\
BUS-7896-0 & 130,395 & 710,935 & 7,883 & 7,883 & 5.6e+09 & 1.0e+09 & 1.5e+10\\
BUS-13160-0 & 140,925 & 1,023,858 & 13,135 & 13,135 & 1.3e+10 & 1.9e+09 & 3.8e+10\\
BUS-18424-0 & 151,455 & 1,338,522 & 18,388 & 18,388 & 2.5e+10 & 2.8e+09 & 7.0e+10\\
BUS-23688-0 & 161,993 & 1,653,437 & 23,643 & 23,643 & 3.9e+10 & 3.8e+09 & 1.1e+11\\
SBT-6255-0 & 163,487 & 809,375 & 6,251 & 6,251 & 5.1e+09 & 1.0e+09 & 1.4e+10\\
SBT-18765-0 & 188,487 & 1,782,614 & 18,745 & 18,745 & 3.3e+10 & 3.5e+09 & 9.6e+10\\
SBT-31275-0 & 213,509 & 2,760,366 & 31,252 & 31,252 & 8.6e+10 & 6.7e+09 & 2.5e+11\\
SBT-43785-0 & 238,529 & 3,741,246 & 43,757 & 43,757 & 1.6e+11 & 1.0e+10 & 4.8e+11\\
SBT-56295-0 & 263,559 & 4,725,174 & 56,269 & 56,269 & 2.7e+11 & 1.5e+10 & 7.8e+11\\
\addlinespace
\bottomrule
\end{tabular}
\end{adjustbox}
\vspace{1em}
\caption{Instance details. Number of variables, constraints, and non-zeros are based on the source-based formulation.}
\label{tbl:instancedetails}
\end{table}

\subsection{Runtime details}
\autoref{tbl:runtime} shows the runtime details for the instances.

\begin{table}[htbp]
\centering
\small
\begin{adjustbox}{max width=\textwidth}
\begin{tabular}{lrrrrrrrrrrr}
\toprule
Instance & \multicolumn{2}{c}{HiGHS} & \multicolumn{2}{c}{MOSEK} & \multicolumn{2}{c}{cuOpt} & \multicolumn{2}{c}{COPT (CPU)} & \multicolumn{2}{c}{COPT (GPU)} & Objective \\
\cmidrule(lr){2-3} \cmidrule(lr){4-5} \cmidrule(lr){6-7} \cmidrule(lr){8-9} \cmidrule(lr){10-11}
 & Tree & Path & Tree & Path & Tree & Path & Tree & Path & Tree & Path & \\
\midrule
\multicolumn{12}{l}{\textbf{Grid}} \\
grid1 & 0.1 & 0.1 & 0.0 & 0.0 & 0.7 & 0.7 & 0.1 & 0.0 & 0.5 & 0.5 & 8.2732e+05 \\
grid2 & 0.1 & 0.1 & 0.0 & 0.0 & 0.9 & 0.7 & 0.1 & 0.1 & 0.6 & 0.4 & 1.7054e+06 \\
grid3 & 0.1 & 0.1 & 0.0 & 0.0 & 0.6 & 0.6 & 0.0 & 0.0 & 0.4 & 0.4 & 1.5246e+06 \\
grid4 & 0.1 & 0.1 & 0.0 & 0.0 & 0.8 & 0.7 & 0.1 & 0.0 & 0.6 & 0.5 & 3.0317e+06 \\
grid5 & 0.1 & 0.1 & 0.0 & 0.0 & 0.8 & 0.7 & 0.1 & 0.1 & 0.6 & 0.5 & 5.0496e+06 \\
grid6 & 0.9 & 0.7 & 0.2 & 0.2 & 1.3 & 1.1 & 0.2 & 0.2 & 1.1 & 0.8 & 1.0401e+07 \\
grid7 & 1.5 & 1.7 & 0.3 & 0.2 & 1.3 & 1.3 & 0.3 & 0.2 & 1.0 & 0.9 & 2.5864e+07 \\
grid8 & 5.3 & 5.6 & 1.4 & 1.2 & 2.5 & 2.2 & 1.5 & 1.6 & 2.6 & 2.3 & 4.1711e+07 \\
grid9 & 14.5 & 16.2 & 3.1 & 2.2 & 3.3 & 2.6 & 3.7 & 3.8 & 3.5 & 2.9 & 8.2653e+07 \\
grid10 & 30.7 & 30.0 & 5.6 & 4.3 & 4.0 & 3.5 & 5.2 & 4.1 & 5.4 & 3.4 & 1.6411e+08 \\
grid11 & 21.4 & 36.0 & 2.9 & 4.6 & 3.3 & 3.1 & 3.6 & 3.4 & 3.0 & 3.4 & 3.2924e+08 \\
grid12 & 16.7 & 36.3 & 3.3 & 4.6 & 2.7 & 2.6 & 2.5 & 3.8 & 2.6 & 3.0 & 5.7713e+08 \\
grid13 & 58.4 & 141.2 & 9.0 & 6.7 & 14.6 & 4.5 & 7.7 & 8.5 & 7.2 & 5.9 & 1.1593e+09 \\
grid14 & 36.8 & 104.8 & 6.5 & 6.8 & 3.6 & 3.9 & 4.0 & 6.8 & 4.6 & 4.7 & 1.8026e+09 \\
grid15 & 101.7 & 268.4 & 10.9 & 16.4 & 16.6 & 8.6 & 11.9 & 19.8 & 9.6 & 9.9 & 3.5935e+09 \\
\addlinespace
\multicolumn{12}{l}{\textbf{Planar}} \\
planar30 & 0.1 & 0.1 & 0.1 & 0.0 & 0.7 & 0.7 & 0.1 & 0.0 & 0.5 & 0.5 & 4.4351e+07 \\
planar50 & 0.2 & 0.2 & 0.1 & 0.1 & 1.1 & 0.9 & 0.1 & 0.1 & 0.9 & 0.5 & 1.2220e+08 \\
planar80 & 1.1 & 1.4 & 0.4 & 0.3 & 1.8 & 1.2 & 0.4 & 0.3 & 1.5 & 0.8 & 1.8244e+08 \\
planar100 & 1.2 & 2.2 & 0.3 & 0.3 & 1.7 & 1.3 & 0.3 & 0.3 & 1.3 & 0.8 & 2.3134e+08 \\
planar150 & 9.1 & 24.3 & 4.9 & 4.6 & 5.4 & 3.6 & 3.8 & 2.9 & 5.2 & 2.7 & 5.4809e+08 \\
planar300 & 5.5 & 11.4 & 1.5 & 1.5 & 2.5 & 2.0 & 1.1 & 1.3 & 2.1 & 1.6 & 6.8998e+08 \\
planar500 & 2.1 & 6.0 & 0.6 & 0.9 & 1.4 & 1.5 & 0.5 & 0.7 & 1.1 & 1.1 & 4.8198e+08 \\
planar800 & 51.9 & 78.6 & 6.1 & 8.0 & 4.5 & 4.5 & 3.8 & 6.7 & 4.9 & 4.9 & 1.1674e+09 \\
planar1000 & 244.1 & 693.0 & 78.8 & 49.0 & 82.2 & 22.2 & 35.1 & 47.2 & 33.6 & 22.8 & 3.4496e+09 \\
planar2500 & t & t & 2071.4 & 908.7 & t & 544.3 & 1122.4 & 1683.4 & 5721.8 & 711.7 & 1.2662e+10 \\
\addlinespace
\multicolumn{12}{l}{\textbf{Transportation}} \\
Austin & 4517.6 & t & 1149.3 & t & 869.8 & t & 558.5 & t & 481.6 & t & 5.4388e+06 \\
Barcelona & 2.5 & 33.2 & 0.7 & 3.5 & 3.3 & 6.4 & 0.8 & 3.6 & 2.8 & 3.6 & 2.4342e+02 \\
BerlinCenter & 45.9 & 665.5 & 14.4 & 21.9 & 34.4 & 22.1 & 11.6 & 15.7 & 11.1 & 12.8 & 2.7529e+07 \\
Birmingham & 2398.2 & t & 564.8 & 1825.0 & 159.5 & 948.1 & 187.5 & 2317.4 & 445.5 & 717.0 & 2.1564e+05 \\
ChicagoRegional & 698.2 & t & 151.8 & 4890.0 & 99.7 & t & 97.3 & 3947.2 & 78.3 & 2819.2 & 5.5132e+06 \\
ChicagoSketch & 1.1 & 44.3 & 0.4 & 5.6 & 1.4 & 5.1 & 0.4 & 3.0 & 1.3 & 3.1 & 6.7059e+06 \\
Philadelphia & 4322.0 & t & 1006.4 & t & 630.1 & t & 425.3 & t & 1535.4 & t & 2.5577e+07 \\
Sydney & 254.5 & t & 75.5 & t & 63.3 & t & 73.6 & t & 66.3 & 6563.0 & 4.9682e+06 \\
Winnipeg & 4.0 & 27.2 & 1.0 & 2.3 & 3.3 & 2.5 & 1.2 & 2.2 & 3.3 & 2.1 & 4.2016e+02 \\
\addlinespace
\multicolumn{12}{l}{\textbf{Intermodal}} \\
BUS-2632-0 & 2.8 & 4.9 & 1.5 & 1.5 & 2.5 & 2.7 & 1.5 & 1.9 & 2.6 & 3.5 & 7.1027e+04 \\
BUS-7896-0 & 8.3 & 7.8 & 2.8 & 2.5 & 3.9 & 4.1 & 2.8 & 2.6 & 3.8 & 3.7 & 2.1059e+05 \\
BUS-13160-0 & 19.1 & 32.3 & 4.1 & 4.7 & 9.1 & 8.2 & 4.5 & 3.7 & 5.4 & 4.6 & 3.4836e+05 \\
BUS-18424-0 & 28.0 & 68.3 & 7.6 & 9.1 & 12.6 & 13.3 & 7.8 & 7.6 & 9.3 & 9.1 & 4.8503e+05 \\
BUS-23688-0 & 68.5 & 34.6 & 10.2 & 9.7 & 11.6 & 14.7 & 8.6 & 8.5 & 9.6 & 9.5 & 6.2488e+05 \\
SBT-6255-0 & 3.2 & 3.7 & 2.2 & 2.0 & 2.9 & 2.9 & 2.0 & 2.6 & 2.5 & 3.2 & 1.6061e+05 \\
SBT-18765-0 & 11.0 & 12.6 & 7.4 & 7.5 & 8.7 & 8.5 & 7.5 & 7.5 & 8.1 & 8.0 & 4.8111e+05 \\
SBT-31275-0 & 22.3 & 24.9 & 16.7 & 21.3 & 22.1 & 22.3 & 20.9 & 16.6 & 21.6 & 16.9 & 8.0132e+05 \\
SBT-43785-0 & 39.1 & 43.0 & 39.2 & 40.7 & 31.1 & 31.2 & 29.8 & 30.1 & 30.5 & 30.6 & 1.1236e+06 \\
SBT-56295-0 & 71.6 & 80.7 & 50.4 & 50.5 & 51.7 & 52.1 & 49.9 & 50.0 & 50.6 & 50.2 & 1.4479e+06 \\
\bottomrule
\end{tabular}
\end{adjustbox}
\vspace{1em}
\caption{Runtime in seconds for the tree- and path-based formulations under each LP backend. A `t' entry marks a run that did not reach the required accuracy within the 7200\,s limit; `-' marks a configuration that produced no result. The objective column reports the best value attained by any run that solved the instance.}
\label{tbl:runtime}
\end{table}

\subsection{Memory details}
\autoref{tbl:memory} shows the peak resident memory for the instances, under both formulations and
all five LP backends. Runs that did not finish are reported at their measured peak rather than at
the machine limit, and flagged as lower bounds; see \autoref{sec:memory}.

\begin{table}[htbp]
\centering
\small
\begin{adjustbox}{max width=\textwidth}
\begin{tabular}{lrrrrrrrrrr}
\toprule
Instance & \multicolumn{2}{c}{HiGHS} & \multicolumn{2}{c}{MOSEK} & \multicolumn{2}{c}{cuOpt} & \multicolumn{2}{c}{COPT (CPU)} & \multicolumn{2}{c}{COPT (GPU)} \\
\cmidrule(lr){2-3} \cmidrule(lr){4-5} \cmidrule(lr){6-7} \cmidrule(lr){8-9} \cmidrule(lr){10-11}
 & Tree & Path & Tree & Path & Tree & Path & Tree & Path & Tree & Path \\
\midrule
\multicolumn{11}{l}{\textbf{Grid}} \\
grid1 & 0.3 & 0.3 & 0.3 & 0.3 & 0.7 & 0.7 & 0.3 & 0.3 & 0.4 & 0.4 \\
grid2 & 0.3 & 0.3 & 0.3 & 0.3 & 0.7 & 0.7 & 0.3 & 0.3 & 0.4 & 0.4 \\
grid3 & 0.3 & 0.3 & 0.3 & 0.3 & 0.7 & 0.7 & 0.3 & 0.3 & 0.4 & 0.4 \\
grid4 & 0.3 & 0.3 & 0.3 & 0.3 & 0.7 & 0.7 & 0.3 & 0.3 & 0.4 & 0.4 \\
grid5 & 0.3 & 0.3 & 0.3 & 0.3 & 0.7 & 0.7 & 0.3 & 0.3 & 0.4 & 0.4 \\
grid6 & 0.3 & 0.3 & 0.3 & 0.3 & 0.7 & 0.7 & 0.3 & 0.3 & 0.5 & 0.5 \\
grid7 & 0.3 & 0.3 & 0.3 & 0.3 & 0.7 & 0.7 & 0.3 & 0.3 & 0.5 & 0.5 \\
grid8 & 0.3 & 0.3 & 0.3 & 0.3 & 0.7 & 0.7 & 0.3 & 0.3 & 0.5 & 0.5 \\
grid9 & 0.3 & 0.3 & 0.3 & 0.3 & 0.7 & 0.7 & 0.3 & 0.3 & 0.5 & 0.5 \\
grid10 & 0.3 & 0.4 & 0.3 & 0.4 & 0.7 & 0.7 & 0.3 & 0.3 & 0.5 & 0.5 \\
grid11 & 0.4 & 0.4 & 0.3 & 0.4 & 0.7 & 0.7 & 0.3 & 0.3 & 0.5 & 0.5 \\
grid12 & 0.4 & 0.4 & 0.3 & 0.3 & 0.7 & 0.7 & 0.3 & 0.3 & 0.5 & 0.5 \\
grid13 & 0.5 & 0.5 & 0.4 & 0.5 & 0.9 & 0.8 & 0.4 & 0.5 & 0.5 & 0.6 \\
grid14 & 0.5 & 0.6 & 0.4 & 0.5 & 0.8 & 0.8 & 0.4 & 0.4 & 0.6 & 0.6 \\
grid15 & 0.6 & 0.9 & 0.5 & 0.8 & 1.1 & 1.0 & 0.4 & 0.7 & 0.6 & 0.9 \\
\addlinespace
\multicolumn{11}{l}{\textbf{Planar}} \\
planar30 & 0.3 & 0.3 & 0.3 & 0.3 & 0.7 & 0.7 & 0.3 & 0.3 & 0.4 & 0.4 \\
planar50 & 0.3 & 0.3 & 0.3 & 0.3 & 0.7 & 0.7 & 0.3 & 0.3 & 0.4 & 0.5 \\
planar80 & 0.3 & 0.3 & 0.3 & 0.3 & 0.7 & 0.7 & 0.3 & 0.3 & 0.5 & 0.5 \\
planar100 & 0.3 & 0.3 & 0.3 & 0.3 & 0.7 & 0.7 & 0.3 & 0.3 & 0.5 & 0.5 \\
planar150 & 0.3 & 0.3 & 0.3 & 0.3 & 0.7 & 0.7 & 0.3 & 0.3 & 0.5 & 0.5 \\
planar300 & 0.3 & 0.3 & 0.3 & 0.3 & 0.7 & 0.7 & 0.3 & 0.3 & 0.5 & 0.5 \\
planar500 & 0.3 & 0.3 & 0.3 & 0.3 & 0.7 & 0.7 & 0.3 & 0.3 & 0.5 & 0.5 \\
planar800 & 0.4 & 0.6 & 0.4 & 0.5 & 0.7 & 0.8 & 0.3 & 0.4 & 0.5 & 0.6 \\
planar1000 & 1.0 & 2.1 & 0.7 & 1.1 & 1.2 & 1.1 & 0.5 & 0.9 & 0.6 & 1.0 \\
planar2500 & 4.0$^\dagger$ & 2.0$^\ast$ & 3.0 & 5.0 & 10.7$^\dagger$ & 4.7 & 1.9 & 4.2 & 2.0 & 4.6 \\
\addlinespace
\multicolumn{11}{l}{\textbf{Transportation}} \\
Austin & 5.4 & 43.3$^\dagger$ & 8.4 & 22.9$^\dagger$ & 10.4 & 7.2$^\ast$ & 4.7 & 16.2$^\dagger$ & 4.7 & 16.0$^\dagger$ \\
Barcelona & 0.3 & 0.4 & 0.4 & 0.3 & 0.7 & 0.7 & 0.3 & 0.3 & 0.5 & 0.5 \\
BerlinCenter & 0.5 & 1.4 & 0.8 & 0.7 & 1.0 & 1.1 & 0.5 & 0.7 & 0.6 & 0.9 \\
Birmingham & 3.4 & 12.2$^\dagger$ & 3.6 & 6.7 & 4.3 & 5.4 & 2.5 & 6.5 & 2.6 & 6.6 \\
ChicagoRegional & 6.1 & 10.9$^\dagger$ & 9.1 & 12.1 & 10.8 & 6.7$^\ast$ & 5.3 & 10.1 & 5.6 & 10.2 \\
ChicagoSketch & 0.4 & 0.5 & 0.4 & 0.5 & 0.8 & 0.9 & 0.3 & 0.6 & 0.5 & 0.7 \\
Philadelphia & 6.7 & 17.5$^\dagger$ & 9.1 & 14.3$^\dagger$ & 11.5 & 7.5$^\ast$ & 5.5 & 9.7$^\dagger$ & 5.5 & 9.4$^\dagger$ \\
Sydney & 9.4 & 21.7$^\dagger$ & 11.2 & 95.8$^\ddagger$ & 13.0 & 10.7$^\ast$ & 9.2 & 21.3$^\dagger$ & 9.4 & 18.6 \\
Winnipeg & 0.3 & 0.4 & 0.4 & 0.3 & 0.7 & 0.7 & 0.3 & 0.3 & 0.5 & 0.5 \\
\addlinespace
\multicolumn{11}{l}{\textbf{Intermodal}} \\
BUS-2632-0 & 0.5 & 0.5 & 0.6 & 0.6 & 0.9 & 0.9 & 0.5 & 0.5 & 0.6 & 0.6 \\
BUS-7896-0 & 0.6 & 0.6 & 0.7 & 0.7 & 1.0 & 1.0 & 0.6 & 0.6 & 0.8 & 0.8 \\
BUS-13160-0 & 0.8 & 0.8 & 0.9 & 0.9 & 1.1 & 1.1 & 0.7 & 0.7 & 0.9 & 0.9 \\
BUS-18424-0 & 0.9 & 0.9 & 1.0 & 1.0 & 1.3 & 1.3 & 0.9 & 0.9 & 1.0 & 1.0 \\
BUS-23688-0 & 1.0 & 1.0 & 1.2 & 1.1 & 1.3 & 1.3 & 1.0 & 1.0 & 1.2 & 1.2 \\
SBT-6255-0 & 0.7 & 0.7 & 0.7 & 0.7 & 1.1 & 1.1 & 0.6 & 0.6 & 0.8 & 0.8 \\
SBT-18765-0 & 1.0 & 1.0 & 1.1 & 1.1 & 1.4 & 1.4 & 1.0 & 1.0 & 1.2 & 1.2 \\
SBT-31275-0 & 1.4 & 1.4 & 1.5 & 1.5 & 1.7 & 1.7 & 1.4 & 1.4 & 1.6 & 1.5 \\
SBT-43785-0 & 1.8 & 1.7 & 1.8 & 1.8 & 2.1 & 2.1 & 1.8 & 1.8 & 1.9 & 1.9 \\
SBT-56295-0 & 2.1 & 2.2 & 2.2 & 2.2 & 2.4 & 2.4 & 2.1 & 2.1 & 2.3 & 2.3 \\
\bottomrule
\end{tabular}
\end{adjustbox}
\vspace{1em}
\caption{Peak resident memory in GB for the tree- and path-based formulations under each LP backend. Every run is measured, including those that did not finish. A $\dagger$ marks a run still running when the 7200\,s limit was reached, whose peak is therefore a lower bound on what it would eventually have needed; $\ast$ marks a run that stopped on its own without reaching the required accuracy, whose peak is a complete measurement of a failed run; and $\ddagger$ marks the one run the kernel terminated for exhausting the machine's 128\,GB. `-' marks a configuration that produced no result. Memory is never imputed to the machine limit, since no run that solved its instance exceeded 18.7\,GB.}
\label{tbl:memory}
\end{table}

\end{document}